\documentclass[a4paper,12pt]{article}
\usepackage{amsmath}
\usepackage[latin1]{inputenc}
\usepackage{amsmath,amsthm,amssymb}
\usepackage{enumerate}
\usepackage[pagewise]{lineno}
\usepackage[colorlinks=true]{hyperref}
\usepackage{cases}

\usepackage{cite}

\hypersetup{urlcolor=blue, citecolor=red}

\renewcommand{\theequation}{\thesection.\arabic{equation}}

\newcommand{\Rmnum}[1]{\expandafter\@slowromancap\romannumeral #1@}
\makeatother
\title{Global solutions to a chemotaxis consumption model involving signal-dependent degenerate diffusion and logistic-type dampening}
\author{Liangchen Wang\footnote{Corresponding author E-mail address: wanglc@cqupt.edu.cn}\\
{\small School of Science, Chongqing University of Posts and Telecommunications,}\\
{\small Chongqing 400065, PR China}}
\date{}
\newtheorem{theorem}{Theorem}[section]

\newtheorem{corollary}{Coroallary}[section]
\newtheorem{lemma}{Lemma}[section]
\newtheorem{definition}{Definition}[section]
\begin{document}
\baselineskip20pt \maketitle
\maketitle
\renewcommand{\theequation}{\arabic{section}.\arabic{equation}}
\catcode`@=11 \@addtoreset{equation}{section} \catcode`@=12

\begin{abstract}
This work considers the Keller-Segel consumption system
\begin{eqnarray*}
\left\{
\begin{array}{llll}
u_t=\Delta (u\phi(v))+au-bu^\gamma,\quad &x\in \Omega,\quad t>0,\\
v_t=\Delta v-uv,\quad &x\in\Omega,\quad t>0
\end{array}
\right.
\end{eqnarray*}
in a smoothly bounded domain $\Omega\subset \mathbb{R}^n,$ $n\geq1$, under no-flux boundary conditions, where the parameters $a,b>0$, $\gamma\geq2$, and the motility function
$\phi$ suitably generalizes the prototype given by $\phi(s)=s^\alpha$ for all $s\geq0$ with $\alpha>0$.

When $\phi$ is appropriately smooth with $\alpha\geq1$, it is shown that if one of the following cases holds:
(i) $\gamma>2$;
(ii) $\gamma=2$, either $n\leq2$ or $n\geq3$ and $b$ is sufficiently large, then for all suitably regular initial data global classical solutions can be constructed.
Whereas when $\phi$ is considered to be with rather mild regularity properties and $\gamma=2$, for arbitrary $b>0$, this system admits at least one global weak solution in case $\alpha>0$. In addition, if $\phi$ is suitably smooth with $\alpha>1$, then the above weak solutions become eventually smooth.
\end{abstract}
\section{Introduction}
\hspace*{\parindent} Chemotaxis is the directed movement of cells as a response to gradients of the concentration of a chemical signal substance, for instance, in the organization of embryonic cell positioning, in tumor growth or also in the governing of immune cell migration.  It is well known to play a crucial role in the interaction of cells with their environment. Since the classical chemotaxis model was proposed by Keller and Segel \cite{KS1,KS2}, many versions of this chemotaxis model have been widely studied (cf. \cite{BBTW,HP,HF}). In particular, a certain type of model where the cell diffusion rate and the chemotactic sensitivity strongly depend on the signal concentration has generated active interest in recent years. As the mathematical core of this model, the diffusion equation reads as follows (see \cite{KS2})
\begin{equation}\label{Equ(1.1)}
u_t=\nabla\cdot(\phi(v)\nabla u-\chi(v)u\nabla v),
\end{equation}
where $u$ and $v$ denote the cell density and the concentration of signals, respectively, and the motility functions $\phi$ and $\chi$ are
linked by a proportionality relation:
$$\chi(v)=(\lambda-1)\phi'(v),$$
where $\lambda$ denotes the ratio of effective body length (i.e., distance between receptors) over the walk length. As announced in \cite{KS2}, the chemotactic motility function $\chi$ may be positive or negative depending on the signs of $\lambda-1$ and $\phi'$.  In addition, this literature also mentioned an interesting special case $\chi(v)=-\phi'(v)$ (i.e. $\lambda=0$), in other words, the distance between receptors is zero and chemotaxis occurs because of an undirected effect on activity due to the presence of a chemical sensed by a single receptor, then (\ref{Equ(1.1)}) can be rewritten as
\begin{equation}\label{Equ(1.01)}
u_t=\Delta(\phi(v)u).
\end{equation}
Such equation has also been proposed in \cite{FTH,LF} to describe processes of stripe pattern formation via so-called ``self-trapping'' mechanisms.

In the recent literature,  there are many works that are highly related to the mathematical analysis of (\ref{Equ(1.01)}). For example, let us mention a Keller-Segel model with signal-dependent motilities, where the signal is produced by cells,
\begin{eqnarray}\label{Equ(1.02)}
\left\{
\begin{array}{llll}
u_t=\Delta (\phi(v)u)+au-bu^\gamma,\quad &x\in \Omega,\quad t>0,\\
v_t=\Delta v-v+u,\quad &x\in\Omega,\quad t>0
\end{array}
\right.
\end{eqnarray}
under no-flux boundary conditions in smoothly bounded domains $\Omega\subset \mathbb{R}^n$ $(n\geq1)$. In the absence of the logistic source (i.e. $a=b=0$), by assuming uniform positive lower and upper bounds on $\phi$, Tao and Winkler \cite{TWM} showed the existences of global bounded classical solutions in two-dimensional case and global weak solutions in higher-dimensional settings, especially, such weak solutions become eventually smooth in three-dimensional counterpart. Note that degeneracy was ruled out due to these assumptions. If the motility function vanishes at large signal concentrations, particularly, by choosing $\phi(s)=s^{-\alpha}$ for all $s>0$ with some $\alpha>0$, global existence of classical solutions with uniform-in-time bounds were discussed in \cite{AYN,JL,FJC,FSN2,FJA,WM,YK}. In addition, global weak solutions with large initial data in lower dimensions $n\leq3$ were constructed (\cite{DKTY}). For simultaneously enhanced diffusion and cross-diffusion, the existence of global weak solutions was discovered in \cite{WN}. While for another special case $\phi(s)=e^{-s}$, it was shown that there exists a critical mass phenomenon in two-dimensional setting (\cite{JW,FJC}). Recently, Fujie and Senba \cite{FSN} proved that global solutions may blow up at infinite time in higher dimensions. Moreover, the existence of global very weak solutions was established in any dimensions (\cite{BLT}). For global solution under another conditions, we refer to \cite{LJ, FJA, WM}.

In the presence of logistic source, global bounded classical solutions and asymptotic behavior were constructed in the quadratic degradation \cite{FJJ,JKW, WW, XL}. Global weak solutions were studied in \cite{DLTW}. For superlinear degradation we refer to \cite{LW,LWL}. Moreover, global existence of classical solutions has also been extensively investigated involving three-component signal-dependent motility models (cf.\cite{JSW,XW,JFZC}).

When chemical is consumed, rather than produced, by cells, the results are far from satisfactory.  More precisely, we consider the following chemotaxis-consumption model involving signal-dependent diffusion
\begin{eqnarray}\label{Equ(1.03)}
\left\{
\begin{array}{llll}
u_t=\Delta (\phi(v)u)+au-bu^\gamma,\quad &x\in \Omega,\quad t>0,\\
v_t=\Delta v-uv,\quad &x\in\Omega,\quad t>0.
\end{array}
\right.
\end{eqnarray}
In the case $a=b=0$, if $\phi\in C^3([0,\infty))$ is strictly positive on $[0,\infty)$, global bounded classical solutions for small $\|v_0\| _{L^\infty(\Omega)}$ and global weak solutions were found in \cite{LZ}. The smallness assumptions of $v_0$ was removed when $n\leq2$ in \cite{LWR}, apart from that, the authors proved that this model has global weak solutions in higher-dimensional cases, especially in three dimension, each of these weak solutions becomes eventually smooth. For some weaker regularity properties on $\phi$, the existence of global very weak solutions was discovered in \cite{LWRR}. When $\phi\in C^1([0,\infty))$ with $\phi(0)=0$ and $\phi'(0)>0$, Winkler \cite{WAA} proved that this system admits global classical solutions for positive times and $u$ approaches a nonconstant profile in the large time limit when $n\leq2$. For enhanced diffusion, global weak solution was shown in \cite{LWP}. In case $\phi(s)=s^{\alpha}$ for all $s\geq0$ with some $\alpha>0$, global very weak solutions in one-and two-dimensional domains were treated in \cite{WAB}. In addition, the existence of global weak solutions was proven for $\alpha\in(0,2]$ in \cite{WAP}. This condition was removed in a global generalized sense (\cite{WZ}). Very recently, for the specific polynomial decay function, i.e. $\phi(s)=s^{-\alpha}$ for all $s>0$ with some $\alpha>0$, both global very weak-strong solutions and global weak-strong solutions were demonstrated in \cite{TWG}. Under the effects of the logistic damping, for some stronger regularity on the initial data, global classical solution was established in \cite{L} if $\gamma>\max\{2,\frac{n+2}{2}\}$. When $\phi\in C^3([0,\infty))$ is strictly positive on $[0,\infty)$, global bounded classical solutions were found in higher-dimensional settings under conditions that either $\gamma>2$ or $\gamma=2$ and $b$ is large, in addition, for all $n\geq3$ and any $b>0$,  global weak solutions were discovered to exist and each of these solutions becomes smooth after some waiting time (\cite{W}).\\
\textbf{Main results.} As discussed in \cite{WAA,WAP}, in the absence of the logistic source, when $\phi(s)=s$, (\ref{Equ(1.03)}) has global classical solutions in case $n\leq2$, and when $\phi(s)=s^\alpha$ with $\alpha\in(0,2]$, global weak solutions were shown for arbitrary dimensions. Hence, the first objective of this manuscript is to explore how far the restrictions of the logistic growth terms such that classical solution exists globally for any dimensions to this model which allows for more general degenerate motilities. A second goal will investigate global weak solutions under rather mild assumptions on $\phi$ and the logistic damping effect. Specifically, we shall consider the initial-boundary value problem
\begin{eqnarray}\label{Equ(1.2)}
\left\{
\begin{array}{llll}
u_t=\Delta (\phi(v)u)+au-bu^\gamma,\quad &x\in \Omega,\quad t>0,\\
v_t=\Delta v-uv,\quad &x\in\Omega,\quad t>0,\\
\frac{\partial (u\phi(v))}{\partial \nu}=\frac{\partial v}{\partial \nu}=0,\quad &x\in\partial\Omega,\,\, t>0,\\
u(x,0)=u_0(x),\,\,\,v(x,0)=v_0(x),\quad &x\in\Omega,\\
\end{array}
\right.
\end{eqnarray}
where $a>0,b>0$, $\gamma\geq2$ and the initial value $u_0$ and $v_0$ satisfy
\begin{equation}\label{Equ(1.3)}
(u_0,v_0)\in [W^{1,\infty}(\Omega)]^2 \text{ with } u_0>0\;\text{and} \,v_0>0\;\text{in} \,\, \overline{\Omega}.
\end{equation}

Our first result will concentrate on global existence of classical solutions.
\begin{theorem}\label{result1.1}
Let $n\geq1$ and $\Omega\subset \mathbb{R}^n$ be a bounded domain with smooth boundary.
Suppose that
\begin{equation}\label{Equ(1.4)}
\phi \in C^1([0, \infty)) \cap C^3((0, \infty)) \text { is such that }  \phi>0 \text { on } (0,\infty),
\end{equation}
and that for some $\alpha\geq1$ we have
\begin{equation}\label{Equ(1.5)}
\liminf _{s \searrow 0} \frac{\phi(s)}{s^\alpha}>0
\end{equation}
and
\begin{equation}\label{Equ(1.6)}
\limsup _{s \searrow 0} \frac{\left|\phi^{\prime}(s)\right|}{s^{\alpha-1}}<\infty.
\end{equation}
Whenever $a>0$, $b>0$ and the initial data $u_0$ and $v_0$ satisfy (\ref{Equ(1.3)}), there exist positive constants $\kappa_1=\kappa_1(n)$, $\kappa_2=\kappa_2(n)$ and $\lambda=\lambda_\phi$ such that if one of the following cases holds:\\
(i) $\gamma>2$;\\
(ii) $\gamma=2$ and $b>\left\{
\begin{array}{llll}
0,\quad &\text{if } n\leq2,\\
\kappa_1(n)\lambda_\phi^{\frac{n+2}{n}}\|v_0\|_{L^\infty(\Omega)}^{\frac{(n+2)\alpha-2}{n}}+\kappa_2(n)\|v_0\|_{L^\infty(\Omega)},\ &\text{if }n\geq3,
\end{array}
\right.$\\
then one can find positive functions
\begin{equation}\label{Equ(1.7)}
\left\{\begin{array}{l}
u \in C^0(\bar{\Omega} \times[0, \infty)) \cap C^{2,1}(\bar{\Omega} \times(0, \infty)),\\
v \in C^0(\bar{\Omega} \times[0, \infty)) \cap C^{2,1}(\bar{\Omega} \times(0, \infty))
\end{array}\right.
\end{equation}
such that $(u, v)$ solves the problem (\ref{Equ(1.2)}) in the classical sense.
\end{theorem}

Theorem \ref{result1.1} shows that stronger logistic degradation and stronger degeneracies are conducive to ensuring that global classical solutions exist to the system (\ref{Equ(1.2)}), especially in higher-dimensional settings.  What can be said about the case of arbitrary $b>0$ and mild assumptions on $\phi$ in the quadratic degradation? In general, it seems that the existence of global classical solution to (\ref{Equ(1.2)}) can not be expected. However, it is possible that for any positive $b$ global weak solutions to the model (\ref{Equ(1.2)}) can be shown to exist regardless of the size of the initial data. Before stating our results, we first give the definition of weak solutions to (\ref{Equ(1.2)}).

\begin{definition}\label{result0.1}
Let $T>0$ and $\gamma=2$. A pair ($u, v$) of nonnegative functions defined on $\Omega \times[0, T)$ is called a weak solution to (\ref{Equ(1.2)}), if
\begin{equation}\label{Equ(1.8)}
\left\{\begin{array}{l}
u \in L^1((0,T); L^1(\Omega)) \quad \text { and } \\
v \in L^\infty(\Omega\times(0, T))  \cap L^{1}\left((0,T); W^{1,1}(\Omega)\right)
\end{array}\right.
\end{equation}
such that
\begin{equation}\label{Equ(1.9)}
u^2 \text{ and } \nabla(u\phi(v))\text{ belong to } L^1((0,T); L^1(\Omega)),
\end{equation}
\begin{equation}\label{Equ(1.10)}
-\int_{0}^{T} \int_{\Omega} u \varphi_{t}-\int_{\Omega} u_0 \varphi (\cdot,0)=-\int_{0}^{T} \int_{\Omega} \nabla(u\phi(v))\cdot\nabla\varphi+a\int_{0}^{T} \int_{\Omega}u\varphi-b\int_{0}^{T} \int_{\Omega}u^2\varphi
\end{equation}
and
\begin{equation}\label{Equ(1.11)}
-\int_{0}^{T} \int_{\Omega} v \varphi_{t}-\int_{\Omega} v_{0} \varphi(\cdot, 0)=-\int_{0}^{T} \int_{\Omega} \nabla v \cdot \nabla \varphi-\int_{0}^{T} \int_{\Omega} u v \varphi
\end{equation}
hold for all $\varphi \in C_{0}^{\infty}(\overline{\Omega }\times[0, T))$.

In particular, if $(u,v)$ is a weak solution of (\ref{Equ(1.2)}) in $\Omega\times(0,T)$ for all $T>0$, then we call $(u,v)$ a global-in-time weak solution to the system (\ref{Equ(1.2)}).
\end{definition}

The existence of global weak solutions will be shown for arbitrary $b>0$ under the restriction of $\alpha\in(0,\frac{1}{2})$.
\begin{theorem}\label{result1.2}
Let $n\geq 1$ and $\Omega\subset \mathbb{R}^n$ be a bounded domain with smooth boundary and let $\gamma=2$, $a>0$ and $b>0$. Suppose that
\begin{equation}\label{Equ(1.12)}
\phi \in C^0([0, \infty)) \cap C^3((0, \infty)) \text { is such that } \phi(0)=0  \text { and }  \phi>0 \text { on } (0,\infty),
\end{equation}
and that with some $\alpha\in(0,\frac{1}{2})$ and $ s_0>0$, apart from (\ref{Equ(1.5)}) and (\ref{Equ(1.6)}),  we have
\begin{equation}\label{Equ(1.13)}
\left(\phi^{\frac{1}{\alpha}}\right)^{\prime \prime}(s) \leq 0 \quad \text { for all } s \in\left(0, s_0\right).
\end{equation}
Then for all $(u_0,v_0)$ fulfilling (\ref{Equ(1.3)}), there exist
\begin{equation}\label{Equ(1.14)}
\left\{\begin{array}{l}
u \in L^\infty((0,\infty); L^1(\Omega))\cap L_{loc}^2([0,\infty); L^2(\Omega)) \quad \text { and } \\
v \in L^\infty(\Omega\times(0, \infty))\cap L^\infty\left((0, \infty) ; W^{1,2}(\Omega)\right)\cap L_{loc}^{2}\left([0, \infty) ; W^{2,2}(\Omega)\right)
\end{array}\right.
\end{equation}
such that $u \geq 0$ and $v>0$ a.e. in $\Omega \times(0, \infty)$, and that $(u,v)$ is a global weak solution of the system (\ref{Equ(1.2)}) in the sense of Definition \ref{result0.1}.
\end{theorem}

Whereas in the case $\alpha\geq\frac{1}{2}$, (\ref{Equ(1.2)}) has global weak solutions for any $b>0$.
\begin{theorem}\label{result1.3}
Let $n\geq 1$ and $\Omega\subset \mathbb{R}^n$ be a bounded domain with smooth boundary and let $\gamma=2$, $a>0$ and $b>0$. Assume that $\alpha\geq\frac{1}{2}$ and that $\phi$ satisfies (\ref{Equ(1.12)}), (\ref{Equ(1.5)}) and (\ref{Equ(1.6)}).
Then whenever $(u_0,v_0)$ complies with (\ref{Equ(1.3)}), we can find functions $u$ and $v$ which satisfy (\ref{Equ(1.14)}) with $u \geq 0$ and $v>0$ a.e. in $\Omega \times(0, \infty)$, and which are such that $(u,v)$ forms a global weak solution of (\ref{Equ(1.2)}) in the sense of Definition \ref{result0.1}.
\end{theorem}

Finally, a slightly more restrictive condition on $\phi$ ensures that the global weak solution gained in Theorem \ref{result1.3} becomes eventually smooth and classical.
\begin{theorem}\label{result1.4}
Let $n\geq 3$ and $\Omega\subset \mathbb{R}^n$ be a bounded domain with smooth boundary and let $\gamma=2$, $a>0$ and $b>0$.
Assume that $\alpha>1$ and that $\phi$ satisfies (\ref{Equ(1.4)})-(\ref{Equ(1.6)}), and that $(u_0,v_0)$ fulfills (\ref{Equ(1.3)}).
Then there exists $T>0$ such that the global weak solution of (\ref{Equ(1.2)}) from Theorem \ref{result1.3} has properties
\begin{equation}\label{Equ(1.15)}
u \in C^{2,1}(\bar{\Omega} \times[T, \infty)) \quad \text { and } \quad v \in C^{2,1}(\bar{\Omega} \times[T, \infty)).
\end{equation}
\end{theorem}

\textbf{Sketch the proof.} In Section 2, we will be concerned with solutions to the approximate system (\ref{Equ(2.1)}) (see Lemma \ref{result2.1}) and prepare some fundamental inequalities (Lemmata \ref{result2.7}, \ref{result2.3} and \ref{result2.4}). Beyond these, in Lemma \ref{result2.5} we derive a differential inequality of $\int_{\Omega}v^{-q+1}_{\varepsilon}(\cdot,t)\left|\nabla v_{\varepsilon}(\cdot,t)\right|^q$ for all $q\geq2$, which will help us to obtain some regularity estimates of $v_\varepsilon$ (Lemma \ref{result2.6}) and be frequently used in the sequel. In Section 3, we will construct functionals in the form $\int_{\Omega} u^p_\varepsilon(\cdot,t)$ for all $p>1$ (Lemma \ref{result3.1}) and $\int_{\Omega} u^p_\varepsilon(\cdot,t)+\int_{\Omega}v^{-q+1}_{\varepsilon}(\cdot,t)\left|\nabla v_{\varepsilon}(\cdot,t)\right|^q$ (Lemma \ref{result3.2}). Based on the differential inequality from Lemma \ref{result3.2}, we will derive $L^p$-bounds of $u_\varepsilon$ in the case $\gamma>2$ (Lemmata \ref{result3.3}-\ref{result3.5}) and in the case $\gamma=2$ and $n\geq3$ if $b$ is sufficiently large (Lemma \ref{result3.6}). Whereas in the case $\gamma=2$ and $n\leq2$, we obtain the boundedness of $u_\varepsilon$ in $L^2$ through a Sobolev inequality and a standard testing procedure (Lemma \ref{result3.7}). In order to overcome the possible degeneracy of diffusion for the first equation in (\ref{Equ(2.1)}), we construct a time-dependent pointwise lower bound for $v_\varepsilon$ (Lemma \ref{result3.8}). Then with the aid of local boundedness criterion (Lemma \ref{result3.9}), we obtain locally bounded global solutions (Lemma \ref{result3.10}). Finally, by the bounds, we derive H\"{o}der estimates and higher order regularity features of $u_\varepsilon$ and $v_\varepsilon$ (Lemmata \ref{result3.11} and \ref{result3.12}), which lead to the conclusion of Theorem \ref{result1.1}. In Section 4, we will devote to the construction of global weak solutions of (\ref{Equ(1.2)}) with arbitrary $b>0$ when $\gamma=2$. Lemma \ref{result4.1} will consider the derivation of estimates about the products $u_\varepsilon v_\varepsilon^\beta$ with some $\beta\geq1$, which allow for compactness arguments and then lead to the pointwise a.e. convergence of $(u_\varepsilon)_{\varepsilon\in(0,1)}$. Because the smaller $\alpha$ will lead to the stronger singularity of cross-diffusion. Hence, the proof of global weak solution will be divided into two cases: The case $\alpha\in(0,\frac{1}{2})$ (Lemmata \ref{result4.2}-\ref{result4.6}) and the case $\alpha\geq\frac{1}{2}$ (Lemmata \ref{result4.7}-\ref{result4.9}). In Section 5, eventual smoothness will be considered: We first construct a functional as follows $\mathcal{E_\varepsilon}(t):= \int_{\Omega}\left(u_\varepsilon(\cdot,t)-\frac{a}{b}-\frac{a}{b}\ln \frac{bu_\varepsilon(\cdot,t)}{a}\right)$
(Lemma \ref{result5.2}). Based on this functional, a uniform lower bound mass of $u_\varepsilon$ (Lemma \ref{result5.3}) will be proved. Thanks to the lower bound for the cell mass, the decay of $v_\varepsilon$ in the $L^\infty$ norm will be shown (Lemma \ref{result5.4}). With the decay of $v_\varepsilon$ at hand, we derive a differential inequality of $\int_{\Omega} u^p_\varepsilon(\cdot,t)+\int_{\Omega}v^{-2p+1}_{\varepsilon}(\cdot,t)\left|\nabla v_{\varepsilon}(\cdot,t)\right|^{2p}$ to establish eventual $L^p$ bounds of $u_\varepsilon$ for all $p>1$ (Lemma \ref{result5.5}). Hence, the higher order regularity of $u_\varepsilon$ and $v_\varepsilon$ can be demonstrated (Lemma \ref{result5.6}).

\section{Preliminaries}
\hspace*{\parindent} In order to appropriately regularize (\ref{Equ(1.2)}), for each $\varepsilon\in(0,1)$, we consider the regularized variant of (\ref{Equ(1.2)}) given by
\begin{eqnarray}\label{Equ(2.1)}
\left\{
\begin{array}{llll}
u_{\varepsilon t}=\Delta (u_{\varepsilon}\phi_\varepsilon(v_\varepsilon))+au_{\varepsilon}-b u^\gamma_{\varepsilon},\quad &x\in \Omega,\quad t>0,\\
v_{\varepsilon t}=\Delta v_\varepsilon-\frac{u_\varepsilon v_\varepsilon}{1+\varepsilon u_\varepsilon},\quad &x\in\Omega,\quad t>0,\\
\frac{\partial u_\varepsilon}{\partial \nu}=\frac{\partial v_\varepsilon}{\partial \nu}=0,\quad &x\in\partial\Omega,\quad t>0,\\
u_\varepsilon(x,0)=u_0(x),\,\,\,v_\varepsilon(x,0)=v_0(x),\quad &x\in\Omega
\end{array}
\right.
\end{eqnarray}
with
\begin{equation}\label{Equ(2.2)}
\phi_\varepsilon(s):=\phi(s)+\varepsilon \quad \text{ for all  } s\geq0 \text{ and }\varepsilon\in(0,1).
\end{equation}
Based on the well-established parabolic theory in \cite{A}, for each $\varepsilon\in(0,1)$, the regularized problem of (\ref{Equ(2.1)}) is globally solvable in the classical sense (see \cite[Lemma 2.2]{WAB}  for details).
\begin{lemma}\label{result2.1}
Let $\Omega\subset \mathbb{R}^{n}(n\geq1)$ be a bounded domain with smooth boundary, and let $a>0$, $b>0$ and $\gamma>1$. Assume that (\ref{Equ(1.3)})
and (\ref{Equ(1.12)}) hold. Then for each $\varepsilon\in(0,1)$, there exist
\begin{eqnarray*}
\left\{
\begin{array}{llll}
u_\varepsilon\in C^0(\overline{\Omega}\times[0,\infty))\cap C^{2,1}(\overline{\Omega}\times(0,\infty)),\\
v_{\varepsilon} \in \bigcap_{q>n} C^{0}\left([0, \infty) ; W^{1, q}(\Omega)\right) \cap C^{2,1}(\bar{\Omega} \times(0, \infty)),
\end{array}
\right.
\end{eqnarray*}
such that $u_\varepsilon>0$ and $v_\varepsilon>0$ in $\overline{\Omega}\times[0,\infty)$, and that $(u_\varepsilon,v_\varepsilon)$ solves (\ref{Equ(2.1)}) classically in $\Omega\times(0,\infty)$.
\end{lemma}

Next we give some elementary estimates.
\begin{lemma}\label{result2.2}
Let $u_0$ and $v_0$ satisfy (\ref{Equ(1.3)}), and suppose that $a>0$, $b>0$ and $\gamma>1$. Then we see that
\begin{equation}\label{Equ(2.3)}
\int_{\Omega} u_\varepsilon(\cdot,t)\leq m_1:=\max \left\{\int_{\Omega} u_{0}(x),\left(\frac{a}{b}\right)^{\frac{1}{\gamma-1}}|\Omega|\right\}
\quad \text{for all  } t>0 \text{ and }\varepsilon\in(0,1),
\end{equation}
\begin{equation}\label{Equ(2.4)}
\int_{t}^{t+1}\int_\Omega u_\varepsilon^\gamma\leq \frac{(a+1)m_1}{b} \quad  \text{for all  } t>0 \text{ and }\varepsilon\in(0,1),
\end{equation}
\begin{equation}\label{Equ(2.5)}
\|v_\varepsilon(\cdot,t)\|_{L^{\infty}(\Omega)}\leq\|v_{0}\|_{L^{\infty}(\Omega)}
\quad  \text{for all  } t>0 \text{ and }\varepsilon\in(0,1)
\end{equation}
and for any $\varepsilon\in(0,1)$,
\begin{equation}\label{Equ(2.6)}
t\mapsto\|v_\varepsilon(\cdot,t)\|_{L^{\infty}(\Omega)}
\quad \text{is nonincreasing in } (0,\infty).
\end{equation}
\end{lemma}
\noindent{\bf{Proof.}} We integrate the first equation in (\ref{Equ(2.1)}) to obtain
\begin{equation}\label{Equ(2.7)}
\begin{aligned}
\frac{d}{dt}\int_{\Omega} u_\varepsilon=a\int_{\Omega}u_\varepsilon-b\int_{\Omega}u_\varepsilon^\gamma \quad  \text{for all  } t>0 \text{ and }\varepsilon\in(0,1).
\end{aligned}
\end{equation}
Thus an ODE comparison argument implies that (\ref{Equ(2.3)}) is valid. And consequently, an integration of (\ref{Equ(2.7)}) shows (\ref{Equ(2.4)}). Both (\ref{Equ(2.5)}) and (\ref{Equ(2.6)}) are a consequence of the maximum principle and the nonnegativity of  $u_\varepsilon$ and $v_\varepsilon$. $\hfill{} \Box$

In the next lemma we will recall an ODE comparison result that will help us to obtain the $L^2$-estimate of $u_\varepsilon$ in case of $\gamma=2$ and $n\leq2$.
\begin{lemma}\label{result2.7} (\cite[Lemma 2.2]{WLP})
Let $T >0$ and $\kappa>1$, and suppose that $c_1>0$ and $c_2 > 0$, and that
$y\in C^0([0, T))\cap C^1((0, T))$ is positive on $[0,T)$ and satisfies
\begin{equation*}
y' (t) + c_1y^\kappa(t) \leq g(t)y(t) \quad \text{for all}\,\,\,t \in (0, T)
\end{equation*}
with some nonnegative function $g(t) \in L^1_{loc}([0, T))$ fulfilling
\begin{equation*}
\int^{t+1}_{t}g(s)ds \leq c_2  \quad \text{for all}\,\,\,t \in (0, T).
\end{equation*}
Then
\begin{equation*}
y(t) \leq \max\left \{y(0)e^{c_2}, (c_1(\kappa-1))^{-\frac{1}{\kappa-1}}e^{c_2}\right\} \quad \text{for all}\,\,\,t \in [0, T).
\end{equation*}
\end{lemma}

The following two crucial lemmas will be frequently used in the sequel.
\begin{lemma}\label{result2.3} (\cite[Lemma 3.4]{WAD})
Let $q\geq2$ and $\psi\in C^2(\overline{\Omega})$ be positive fulfilling $\frac{\partial\psi}{\partial\nu}=0$ on $\partial\Omega$. Then we have
\begin{equation*}
\int_{\Omega}\psi^{-q-1}\left|\nabla \psi\right|^{q+2}
\leq(q+\sqrt{n})^2\int_{\Omega} \psi^{-q+3}\left|\nabla \psi\right|^{q-2}\left|D^2 \ln \psi\right|^2
\end{equation*}
and
\begin{equation*}
\int_{\Omega}\psi^{-q+1}\left|\nabla \psi\right|^{q-2}\left|D^2 \psi\right|^2
\leq(q+\sqrt{n}+1)^2\int_{\Omega} \psi^{-q+3}\left|\nabla \psi\right|^{q-2}\left|D^2 \ln \psi\right|^2.
\end{equation*}
\end{lemma}

\begin{lemma}\label{result2.4} (\cite[Lemma 3.5]{WAD})
Let $q\geq2$ and $\eta>0$. There is $C=C(q,\eta)>0$ such that for every positive $\psi\in C^2(\overline{\Omega})$ with $\frac{\partial\psi}{\partial\nu}=0$ on $\partial\Omega$ satisfies
\begin{equation*}
\int_{\partial \Omega} \psi^{-q+1}\left|\nabla \psi\right|^{q-2}\frac{\partial\left|\nabla \psi\right|^2}{\partial \nu}\leq\eta\int_{\Omega}\psi^{-q-1}\left|\nabla \psi\right|^{q+2}+\eta\int_{\Omega}\psi^{-q+1}\left|\nabla \psi\right|^{q-2}\left|D^2 \psi\right|^2+C\int_{\Omega}\psi.
\end{equation*}
\end{lemma}

Below, we establish a differential inequality of $\int_{\Omega} v_{\varepsilon}^{-q+1}(\cdot,t)\left|\nabla v_{\varepsilon}(\cdot,t)\right|^q$ for all $q\geq2$ and $\varepsilon\in(0,1)$, which will be referred to several times throughout the sequel and which is an extension of \cite[Lemma 3.3]{WAD} where the reaction part in the second equation is not regularized, whereas that is regularized in our model. But the idea of our proof is essentially inspired by \cite[Lemma 3.3]{WAD} and we provide necessary details of the proof below for clarity.
\begin{lemma}\label{result2.5}
If (\ref{Equ(1.3)}) holds, then for all $q\geq2$ we obtain
\begin{equation}\label{Equ(2.8)}
\begin{aligned}
\frac{d}{d t} \int_{\Omega} v_{\varepsilon}^{-q+1}\left|\nabla v_{\varepsilon}\right|^q+q(q-1)\int_{\Omega} v_{\varepsilon}^{-q+3}&\left|\nabla v_{\varepsilon}\right|^{q-2}\left|D^2 \ln v_{\varepsilon}\right|^2\leq\frac{q}{2} \int_{\partial \Omega} v_{\varepsilon}^{-q+1}\left|\nabla v_{\varepsilon}\right|^{q-2} \frac{\partial\left|\nabla v_{\varepsilon}\right|^2}{\partial \nu} \\
&+q(q-2+\sqrt{n}) \int_{\Omega} u_{\varepsilon} v_{\varepsilon}^{-q+2}\left|\nabla v_{\varepsilon}\right|^{q-2}\left|D^2 v_{\varepsilon}\right|
\end{aligned}
\end{equation}
for all $t>0$ and $\varepsilon\in(0,1)$.
\end{lemma}
\noindent{\bf{Proof.}} Integrating by parts to the second equation in (\ref{Equ(2.1)}) and using $2\nabla v_{\varepsilon}\cdot\nabla\Delta v_{\varepsilon}=\Delta |\nabla v_{\varepsilon}|^2-2|D^2v_{\varepsilon}|^2$, we know that
\begin{equation}\label{Equ(2.9)}
\begin{aligned}
\frac{d}{d t} \int_{\Omega} v_{\varepsilon}^{-q+1}\left|\nabla v_{\varepsilon}\right|^q=& q \int_{\Omega} v_{\varepsilon}^{-q+1}\left|\nabla v_{\varepsilon}\right|^{q-2} \nabla v_{\varepsilon} \cdot\nabla\left\{ \Delta v_{\varepsilon}-\frac{u_{\varepsilon} v_{\varepsilon}}{1+\varepsilon u_\varepsilon}\right\} \\
&-(q-1) \int_{\Omega} v_{\varepsilon}^{-q}\left|\nabla v_{\varepsilon}\right|^q\left\{\Delta v_{\varepsilon}-\frac{u_{\varepsilon} v_{\varepsilon}}{1+\varepsilon u_\varepsilon}\right\}\\
=& \frac{q}{2} \int_{\Omega} v_{\varepsilon}^{-q+1}\left|\nabla v_{\varepsilon}\right|^{q-2} \Delta\left|\nabla v_{\varepsilon}\right|^2-q \int_{\Omega} v_{\varepsilon}^{-q+1}\left|\nabla v_{\varepsilon}\right|^{q-2}\left|D^2 v_{\varepsilon}\right|^2 \\
&-q \int_{\Omega} v_{\varepsilon}^{-q+1}\left|\nabla v_{\varepsilon}\right|^{q-2} \nabla v_{\varepsilon} \cdot \nabla\left(\frac{u_{\varepsilon} v_{\varepsilon}}{1+\varepsilon u_\varepsilon}\right) \\
&-(q-1) \int_{\Omega} v_{\varepsilon}^{-q}\left|\nabla v_{\varepsilon}\right|^q \Delta v_{\varepsilon}+(q-1) \int_{\Omega} \frac{u_{\varepsilon} v_{\varepsilon}^{-q+1}\left|\nabla v_{\varepsilon}\right|^q}{1+\varepsilon u_\varepsilon}\\
=&q(q-1) \int_{\Omega} v_{\varepsilon}^{-q}\left|\nabla v_{\varepsilon}\right|^{q-2} \nabla v_{\varepsilon} \cdot \nabla\left|\nabla v_{\varepsilon}\right|^2-q \int_{\Omega} v_{\varepsilon}^{-q+1}\left|\nabla v_{\varepsilon}\right|^{q-2}\left|D^2 v_{\varepsilon}\right|^2\\
&-\frac{q(q-2)}{4} \int_{\Omega} v_{\varepsilon}^{-q+1}\left|\nabla v_{\varepsilon}\right|^{q-4}\left|\nabla| \nabla v_{\varepsilon}|^2\right|^2
-q(q-1) \int_{\Omega} v_{\varepsilon}^{-q-1}\left|\nabla v_{\varepsilon}\right|^{q+2} \\
&+\frac{q}{2} \int_{\partial \Omega} v_{\varepsilon}^{-q+1}\left|\nabla v_{\varepsilon}\right|^{q-2} \cdot \frac{\partial\left|\nabla v_{\varepsilon}\right|^2}{\partial \nu}\\
&+\frac{q(q-2)}{2} \int_{\Omega} \frac{u_{\varepsilon} v_{\varepsilon}^{-q+2}}{1+\varepsilon u_\varepsilon} \left|\nabla v_{\varepsilon}\right|^{q-4} \nabla v_{\varepsilon} \cdot\nabla| \nabla v_{\varepsilon}|^2 \\
&+q \int_{\Omega} \frac{u_{\varepsilon} v_{\varepsilon}^{-q+2}}{1+\varepsilon u_\varepsilon}\left|\nabla v_{\varepsilon}\right|^{q-2} \Delta v_{\varepsilon}-(q-1)^2 \int_{\Omega} \frac{u_{\varepsilon} v_{\varepsilon}^{-q+1}}{1+\varepsilon u_\varepsilon} \left|\nabla v_{\varepsilon}\right|^q\\
\end{aligned}
\end{equation}
for all $t>0$ and $\varepsilon\in(0,1)$. For the first four terms on the right of (\ref{Equ(2.9)}), we make use of the pointwise identity (\cite[Lemma 3.2]{LM2} or \cite[Lemma 3.2]{WAD})
$$
\left|D^2 \ln \varphi\right|^2=-\frac{1}{\varphi^3} \nabla \varphi \cdot \nabla|\nabla \varphi|^2+\frac{1}{\varphi^2}\left|D^2 \varphi\right|^2+\frac{1}{\varphi^4}|\nabla \varphi|^4 \quad \text { for all positive } \varphi \in C^2(\overline{\Omega})
$$
and $\nabla|\nabla v_{\varepsilon}|^2=2D^2 v_{\varepsilon} \cdot \nabla v_{\varepsilon}$ to obtain
\begin{equation}\label{Equ(2.10)}
\begin{split}
&q(q-1) \int_{\Omega} v_{\varepsilon}^{-q}\left|\nabla v_{\varepsilon}\right|^{q-2} \nabla v_{\varepsilon} \cdot \nabla\left|\nabla v_{\varepsilon}\right|^2-\frac{q(q-2)}{4} \int_{\Omega} v_{\varepsilon}^{-q+1}\left|\nabla v_{\varepsilon}\right|^{q-4}\left|\nabla| \nabla v_{\varepsilon}|^2\right|^2 \\
&-q \int_{\Omega} v_{\varepsilon}^{-q+1}\left|\nabla v_{\varepsilon}\right|^{q-2}\left|D^2 v_{\varepsilon}\right|^2-q(q-1) \int_{\Omega} v_{\varepsilon}^{-q-1}\left|\nabla v_{\varepsilon}\right|^{q+2}\\
=&-q(q-1) \int_{\Omega} v_{\varepsilon}^{-q+3}\left|\nabla v_{\varepsilon}\right|^{q-2}\Big(-\frac{1}{v_{\varepsilon}^3}\nabla v_{\varepsilon} \cdot\nabla|\nabla v_{\varepsilon}|^2+\frac{1}{v_{\varepsilon}^2}|D^2 v_{\varepsilon}|^2+\frac{1}{v_{\varepsilon}^4}|\nabla v_{\varepsilon}|^4\Big)\\
=&-q(q-1) \int_{\Omega} v_{\varepsilon}^{-q+3}\left|\nabla v_{\varepsilon}\right|^{q-2}\left|D^2 \ln v_{\varepsilon}\right|^2
\end{split}
\end{equation}
for all $t>0$ and $\varepsilon\in(0,1)$. The sixth and the seventh summands on the right of  (\ref{Equ(2.9)}) can be estimated according to
\begin{equation}\label{Equ(2.11)}
\begin{split}
&\frac{q(q-2)}{2} \int_{\Omega} \frac{u_{\varepsilon} v_{\varepsilon}^{-q+2}}{1+\varepsilon u_\varepsilon} \left|\nabla v_{\varepsilon}\right|^{q-4} \nabla v_{\varepsilon} \cdot\nabla| \nabla v_{\varepsilon}|^2+q \int_{\Omega} \frac{u_{\varepsilon} v_{\varepsilon}^{-q+2}}{1+\varepsilon u_\varepsilon}\left|\nabla v_{\varepsilon}\right|^{q-2} \Delta v_{\varepsilon}\\
\leq& q(q-2+\sqrt{n}) \int_{\Omega} u_{\varepsilon} v_{\varepsilon}^{-q+2}\left|\nabla v_{\varepsilon}\right|^{q-2}\left|D^2 v_{\varepsilon}\right|
\end{split}
\end{equation}
for all $t>0$ and $\varepsilon\in(0,1)$, because of $|\Delta v_{\varepsilon}|^2\leq n|D^2v_{\varepsilon}|^2$. Inserting (\ref{Equ(2.10)}) and (\ref{Equ(2.11)}) into (\ref{Equ(2.9)}) we obtain (\ref{Equ(2.8)}). $\hfill{} \Box$

Next, some regularity information on $v_\varepsilon$ can be established.
\begin{lemma}\label{result2.6}
If (\ref{Equ(1.3)}) holds and $\gamma\geq2$, then there exists $C>0$ such that
\begin{align}
\int_\Omega\frac{|\nabla v_\varepsilon(\cdot,t)|^2}{v_\varepsilon(\cdot,t)} \leq C \quad\text{ for all  } t>0 \text{ and }\varepsilon\in(0,1),\label{Equ(2.12)}\\
\int_{t}^{t+1} \int_\Omega \frac{|D^2 v_{\varepsilon}|^2}{v_{\varepsilon}} \leq C \quad\text{ for all } t>0 \text{ and }\varepsilon\in(0,1),\label{Equ(2.13)}\\
\int_{t}^{t+1} \int_\Omega \frac{|\nabla v_\varepsilon|^4}{v^3_\varepsilon} \leq C \quad\text{ for all } t>0 \text{ and }\varepsilon\in(0,1) \text{ and }\label{Equ(2.15)}\\
\int_{t}^{t+1} \int_\Omega v^2_{\varepsilon t}\leq C \quad\text{ for all } t>0 \text{ and }\varepsilon\in(0,1).\label{Equ(2.14)}
\end{align}
Moreover, for any $T>0$, there exists $C(T)>0$ such that
\begin{align}
\int_\Omega \ln \frac{1}{v_\varepsilon(\cdot,t)}\leq C(T)\quad\text{ for all  } t\in(0,T) \text{ and }\varepsilon\in(0,1) \text{ and }\label{Equ(2.16)}\\
\int_{0}^{T} \int_\Omega \frac{|\nabla v_\varepsilon|^2}{v^2_\varepsilon}\leq C(T) \quad\text{ for all } t\in(0,T) \text{ and }\varepsilon\in(0,1).\label{Equ(2.17)}
\end{align}
\end{lemma}
\noindent{\bf{Proof.}} Lemma \ref{result2.5} can be applied to $q=2$ to show that
\begin{equation}\label{Equ(2.18)}
\begin{aligned}
\frac{d}{d t} \int_{\Omega}\frac{|\nabla v_{\varepsilon}|^2}{v_{\varepsilon}}+2\int_{\Omega} v_{\varepsilon}|D^2 \ln v_{\varepsilon}|^2\leq&\int_{\partial \Omega} \frac{1}{v_{\varepsilon}}\frac{\partial\left|\nabla v_{\varepsilon}\right|^2}{\partial \nu}+2\sqrt{n} \int_{\Omega} u_{\varepsilon}\left|D^2 v_{\varepsilon}\right|
\end{aligned}
\end{equation}
for all $t>0$ and $\varepsilon\in(0,1)$. An application of Lemma \ref{result2.4} and (\ref{Equ(2.5)}) shows the existence of $c_1>0$ such that
\begin{equation}\label{Equ(2.19)}
\begin{aligned}
\int_{\partial \Omega} \frac{1}{v_{\varepsilon}}\frac{\partial\left|\nabla v_{\varepsilon}\right|^2}{\partial \nu}\leq&\frac{1}{4(2+\sqrt{n})^2}\int_{\Omega}\frac{\left|\nabla v_{\varepsilon}\right|^{4}}{v^3_{\varepsilon}}+\frac{1}{4(3+\sqrt{n})^2}\int_{\Omega}\frac{|D^2 v_{\varepsilon}|^2}{v_{\varepsilon}}+c_1
\end{aligned}
\end{equation}
for all $t>0$ and $\varepsilon\in(0,1)$.
We utilize (\ref{Equ(2.5)}) and Young's inequality to obtain $c_2>0$ fulfilling
\begin{equation}\label{Equ(2.20)}
\begin{aligned}
2\sqrt{n} \int_{\Omega} u_{\varepsilon}\left|D^2 v_{\varepsilon}\right|\leq \frac{1}{4(3+\sqrt{n})^2} \int_{\Omega}\frac{|D^2 v_{\varepsilon}|^2}{v_{\varepsilon}}+4n(3+\sqrt{n})^2\|v_0\|_{L^\infty(\Omega)} \int_{\Omega}u^2_{\varepsilon}
\end{aligned}
\end{equation}
and
\begin{equation}\label{Equ(2.21)}
\begin{aligned}
\int_{\Omega}\frac{|\nabla v_{\varepsilon}|^2}{v_{\varepsilon}}\leq \frac{1}{4(2+\sqrt{n})^2} \int_{\Omega}\frac{\left|\nabla v_{\varepsilon}\right|^{4}}{v^3_{\varepsilon}}+c_2\int_{\Omega}v_\varepsilon \leq \frac{1}{4(2+\sqrt{n})^2} \int_{\Omega}\frac{\left|\nabla v_{\varepsilon}\right|^{4}}{v^3_{\varepsilon}}+c_2|\Omega|\|v_0\|_{L^\infty(\Omega)}
\end{aligned}
\end{equation}
for all $t>0$ and $\varepsilon\in(0,1)$.
It follows from Lemma \ref{result2.3} that
\begin{equation}\label{Equ(2.22)}
\begin{aligned}
\frac{1}{(2+\sqrt{n})^2}\int_{\Omega}\frac{\left|\nabla v_{\varepsilon}\right|^{4}}{v^3_{\varepsilon}}\leq\int_{\Omega} v_{\varepsilon}|D^2 \ln v_{\varepsilon}|^2 \text{ and }
\frac{1}{(3+\sqrt{n})^2}\int_{\Omega}\frac{|D^2 v_{\varepsilon}|^2}{v_{\varepsilon}}\leq\int_{\Omega} v_{\varepsilon}|D^2 \ln v_{\varepsilon}|^2
\end{aligned}
\end{equation}
for all $t>0$ and $\varepsilon\in(0,1)$. Collecting (\ref{Equ(2.18)})-(\ref{Equ(2.22)}), we can find some $c_3>0$ and $c_4>0$ such that
\begin{equation}\label{Equ(2.23)}
\begin{aligned}
\frac{d}{d t} \int_{\Omega}\frac{|\nabla v_{\varepsilon}|^2}{v_{\varepsilon}}+\int_{\Omega}\frac{|\nabla v_{\varepsilon}|^2}{v_{\varepsilon}}+c_3\int_{\Omega}\frac{\left|\nabla v_{\varepsilon}\right|^{4}}{v^3_{\varepsilon}}+c_3\int_{\Omega}\frac{|D^2 v_{\varepsilon}|^2}{v_{\varepsilon}}
\leq c_4\int_{\Omega} u^2_{\varepsilon}+c_4
\end{aligned}
\end{equation}
for all $t>0$ and $\varepsilon\in(0,1)$. This along with \cite[Lemma 3.4]{SSW} and (\ref{Equ(2.4)}) implies (\ref{Equ(2.12)}). Then integrating (\ref{Equ(2.23)}) over $(t,t+1)$, we obtain (\ref{Equ(2.13)}) and (\ref{Equ(2.15)}). The second equation in (\ref{Equ(2.1)}) together with (\ref{Equ(2.4)}), (\ref{Equ(2.5)}) and (\ref{Equ(2.13)}) yields  (\ref{Equ(2.14)}). The last two assertions can be proved in \cite[Lemma 4.1]{WAB}. $\hfill{} \Box$

\section{Proof of Theorem 1.1}

\hspace*{\parindent} In this section, our goal is to obtain global classical solutions of (\ref{Equ(1.2)}).  To this end, we will establish certain integral inequalities of (\ref{Equ(2.1)}), which are of crucial importance to achieve the $L^p$ regularity argument concerning $u_\varepsilon$. The ideas used in this section mainly come from \cite{BBTW,LW1,WAA,WAD}. We first use a standard testing procedure for the first equation in (\ref{Equ(2.1)}).

\begin{lemma}\label{result3.1}
Suppose that (\ref{Equ(1.5)}) and (\ref{Equ(1.6)}) are valid. Then there exist $C_1>0$ and $C_2>0$ such that whenever (\ref{Equ(1.3)}) holds,
for all $p>1$ we have
\begin{equation}\label{Equ(3.1)}
\begin{aligned}
\frac{d}{dt}\int_{\Omega} u_\varepsilon^{p}+\frac{C_1p(p-1)}{2}\int_\Omega u_\varepsilon^{p-2}v_\varepsilon^\alpha|\nabla u_\varepsilon|^2\leq& \frac{C_2p(p-1)}{2}\int_{\Omega} u_\varepsilon^{p}v_\varepsilon^{\alpha-2}|\nabla v_\varepsilon|^2\\
&+ap\int_\Omega u_\varepsilon^{p}-bp\int_{\Omega} u_\varepsilon^{p+\gamma-1}
\end{aligned}
\end{equation}
for all $t>0$ and $\varepsilon\in(0,1)$.
\end{lemma}
\noindent{\bf{Proof.}} Thanks to (\ref{Equ(1.5)}), (\ref{Equ(1.6)}), (\ref{Equ(2.2)}) and (\ref{Equ(2.5)}), we infer that there exist certain positive constants $c_1$ and $c_2$ such that
\begin{equation}\label{Equ(3.2)}
\begin{aligned}
c_1 v^\alpha_\varepsilon\leq\phi_\varepsilon(v_\varepsilon)\,\, \text{ and } |\phi_\varepsilon'(v_\varepsilon)|\leq c_{2} v^{\alpha-1}_\varepsilon \quad \text{ in } \Omega\times(0,\infty) \,\,\text{  for all  } \varepsilon\in(0,1).
\end{aligned}
\end{equation}
We test the first equation in (\ref{Equ(2.1)}) by $u_\varepsilon^{p-1}$ and use Young's inequality to see that
\begin{equation}\label{Equ(3.4)}
\begin{aligned}
\frac{d}{dt}\int_{\Omega} u_\varepsilon^{p}=&-p(p-1)\int_\Omega u_\varepsilon^{p-2}\phi_\varepsilon(v_\varepsilon)|\nabla u_\varepsilon|^2-p(p-1)\int_{\Omega} u_\varepsilon^{p-1}\phi_\varepsilon'(v_\varepsilon)\nabla u_\varepsilon\cdot\nabla u_\varepsilon\\
&+ap\int_\Omega u_\varepsilon^{p}-bp\int_{\Omega} u_\varepsilon^{p+\gamma-1}\\
\leq&-\frac{p(p-1)}{2}\int_\Omega u_\varepsilon^{p-2}\phi_\varepsilon(v_\varepsilon)|\nabla u_\varepsilon|^2+\frac{p(p-1)}{2}\int_{\Omega} \frac{|\phi_\varepsilon'(v_\varepsilon)|^2}{\phi_\varepsilon(v_\varepsilon)}u_\varepsilon^p|\nabla v_\varepsilon|^2\\
&+ap\int_\Omega u_\varepsilon^{p}-bp\int_{\Omega} u_\varepsilon^{p+\gamma-1}
\end{aligned}
\end{equation}
for all $t>0$ and $\varepsilon\in(0,1)$. We thus easily obtain (\ref{Equ(3.1)}) if we let $C_1:=c_1$ and $C_2:=\frac{c_2^2}{c_1}$. $\hfill{} \Box$

By suitably using Lemmata \ref{result2.5} and \ref{result3.1}, we will construct the following quasi-energy functional of the form $\int_{\Omega}u_\varepsilon^p(\cdot,t)+\int_{\Omega} v_{\varepsilon}^{-q+1}(\cdot,t)\left|\nabla v_{\varepsilon}(\cdot,t)\right|^q$ with $p>1$ and $q\geq2$, which will be used in the proof of Lemmata \ref{result3.5}, \ref{result3.6} and \ref{result5.5}.
\begin{lemma}\label{result3.2}
Assume that (\ref{Equ(1.3)}), (\ref{Equ(1.5)}) and (\ref{Equ(1.6)}) are valid with some $\alpha\geq1$, and that $\gamma>1$. Then for all $p>1$ and $q\geq2$ there exists a constant $C_1>0$ such that
\begin{equation}\label{Equ(3.5)}
\begin{aligned}
\frac{d}{dt}&\left(\int_{\Omega}u_\varepsilon^p+\int_{\Omega} v_{\varepsilon}^{-q+1}\left|\nabla v_{\varepsilon}\right|^q\right)+\int_{\Omega}u_\varepsilon^p+\int_{\Omega} v_{\varepsilon}^{-q+1}\left|\nabla v_{\varepsilon}\right|^q\\
\leq& \frac{p(p-1)}{2}\left(\frac{q(q-1)}{2p(p-1)(q+\sqrt{n})^2}\right)^{-\frac{2}{q}}C_2^{\frac{q+2}{q}}
\|v_\varepsilon(\cdot,t)\|_{L^\infty(\Omega)}^\frac{\alpha(q+2)-2}{q}\int_{\Omega}u_\varepsilon^\frac{p(q+2)}{q}\\
&+\frac{8^\frac{q}{2}q(q+\sqrt{n}+1)^q(q-2+\sqrt{n})^\frac{q+2}{2}\|v_\varepsilon(\cdot,t)\|_{L^\infty(\Omega)}}{(q-1)^\frac{q}{2}}
\int_{\Omega}u_\varepsilon^{\frac{q+2}{2}}\\
&-\frac{bp}{2}\int_{\Omega} u_\varepsilon^{p+\gamma-1}+C_1
\end{aligned}
\end{equation}
for all $t>0$ and $\varepsilon\in(0,1)$, where $C_2$ is provided by (\ref{Equ(3.1)}).
\end{lemma}
\noindent{\bf{Proof.}} We only need to take an evident linear combination of  (\ref{Equ(3.1)}) and (\ref{Equ(2.8)}) to obtain
\begin{equation}\label{Equ(3.6)}
\begin{aligned}
\frac{d}{dt}&\left(\int_{\Omega}u_\varepsilon^p+\int_{\Omega} v_{\varepsilon}^{-q+1}\left|\nabla v_{\varepsilon}\right|^q\right)+\int_{\Omega}u_\varepsilon^p+\int_{\Omega} v_{\varepsilon}^{-q+1}\left|\nabla v_{\varepsilon}\right|^q\\
&+q(q-1)\int_{\Omega} v_{\varepsilon}^{-q+3}\left|\nabla v_{\varepsilon}\right|^{q-2}\left|D^2 \ln v_{\varepsilon}\right|^2\\
\leq& \frac{C_2p(p-1)}{2}\int_{\Omega} u_\varepsilon^{p}v_\varepsilon^{\alpha-2}|\nabla v_\varepsilon|^2+
(ap+1)\int_\Omega u_\varepsilon^{p}-bp\int_{\Omega} u_\varepsilon^{p+\gamma-1}\\
&+q(q-2+\sqrt{n}) \int_{\Omega} u_{\varepsilon} v_{\varepsilon}^{-q+2}\left|\nabla v_{\varepsilon}\right|^{q-2}\left|D^2 v_{\varepsilon}\right|\\
&+\frac{q}{2} \int_{\partial \Omega} v_{\varepsilon}^{-q+1}\left|\nabla v_{\varepsilon}\right|^{q-2} \cdot \frac{\partial\left|\nabla v_{\varepsilon}\right|^2}{\partial \nu}+\int_{\Omega} v_{\varepsilon}^{-q+1}\left|\nabla v_{\varepsilon}\right|^q \\
\end{aligned}
\end{equation}
for all $t>0$ and $\varepsilon\in(0,1)$ with $C_2$ taken from Lemma \ref{result3.1}. We employ Young's inequality with any $\eta_1>0$ and Lemma \ref{result2.3} to estimate
\begin{equation}\label{Equ(3.7)}
\begin{aligned}
\frac{C_2p(p-1)}{2}\int_{\Omega} u_\varepsilon^{p}v_\varepsilon^{\alpha-2}|\nabla v_\varepsilon|^2\leq&\frac{p(p-1)}{2}\eta_1^\frac{q+2}{2}\int_{\Omega}v_{\varepsilon}^{-q-1}\left|\nabla v_{\varepsilon}\right|^{q+2}\\
&+\frac{p(p-1)}{2}\eta_1^{-\frac{q+2}{q}}C_2^{\frac{q+2}{q}}\int_{\Omega}u_\varepsilon^\frac{p(q+2)}{q}v_\varepsilon^\frac{\alpha(q+2)-2}{q}\\
\leq&\frac{p(p-1)(q+\sqrt{n})^2}{2}\eta_1^\frac{q+2}{2}\int_{\Omega} v_{\varepsilon}^{-q+3}\left|\nabla v_{\varepsilon}\right|^{q-2}\left|D^2 \ln v_{\varepsilon}\right|^2\\
&+\frac{p(p-1)}{2}\eta_1^{-\frac{q+2}{q}}C_2^{\frac{q+2}{q}}\|v_\varepsilon(\cdot,t)\|_{L^\infty(\Omega)}^\frac{\alpha(q+2)-2}{q}\int_{\Omega}u_\varepsilon^\frac{p(q+2)}{q}
\end{aligned}
\end{equation}
for all $t>0$ and $\varepsilon\in(0,1)$. We pick $\eta_1=\left(\frac{q(q-1)}{2p(p-1)(q+\sqrt{n})^2}\right)^{\frac{2}{q+2}}$ and use (\ref{Equ(3.7)})  to see that
\begin{equation}\label{Equ(3.8)}
\begin{aligned}
\frac{C_2p(p-1)}{2}&\int_{\Omega} u_\varepsilon^{p}v_\varepsilon^{\alpha-2}|\nabla v_\varepsilon|^2\leq\frac{q(q-1)}{4}\int_{\Omega} v_{\varepsilon}^{-q+3}\left|\nabla v_{\varepsilon}\right|^{q-2}\left|D^2 \ln v_{\varepsilon}\right|^2\\
&+\frac{p(p-1)}{2}\left(\frac{q(q-1)}{2p(p-1)(q+\sqrt{n})^2}\right)^{-\frac{2}{q}}C_2^{\frac{q+2}{q}}
\|v_\varepsilon(\cdot,t)\|_{L^\infty(\Omega)}^\frac{\alpha(q+2)-2}{q}\int_{\Omega}u_\varepsilon^\frac{p(q+2)}{q}
\end{aligned}
\end{equation}
for all $t>0$ and $\varepsilon\in(0,1)$. An application of Young's inequality with any $\eta_2, \eta_3>0$  and Lemma \ref{result2.3} to the fourth term on the right hand of (\ref{Equ(3.6)}) entails that
\begin{equation}\label{Equ(3.9)}
\begin{aligned}
&q(q-2+\sqrt{n}) \int_{\Omega} u_{\varepsilon} v_{\varepsilon}^{-q+2}\left|\nabla v_{\varepsilon}\right|^{q-2}\left|D^2 v_{\varepsilon}\right|\\
\leq& \eta_2\int_{\Omega}v_{\varepsilon}^{-q+1}\left|\nabla v_{\varepsilon}\right|^{q-2}\left|D^2 v_{\varepsilon}\right|^2+\eta_2^{-1}q^2(q-2+\sqrt{n})^2 \int_{\Omega}u_\varepsilon^{2}v_{\varepsilon}^{-q+3}\left|\nabla v_{\varepsilon}\right|^{q-2}\\
\leq&\eta_2(q+\sqrt{n}+1)^2\int_{\Omega} v_{\varepsilon}^{-q+3}\left|\nabla v_{\varepsilon}\right|^{q-2}\left|D^2 \ln v_{\varepsilon}\right|^2+\eta_2^{-1}q^2(q-2+\sqrt{n})^2\eta_3^\frac{q+2}{q-2}\int_{\Omega}v_{\varepsilon}^{-q-1}\left|\nabla v_{\varepsilon}\right|^{q+2}\\
&+\eta_2^{-1}q^2(q-2+\sqrt{n})^2\eta_3^{-\frac{q+2}{4}} \int_{\Omega}u_\varepsilon^{\frac{q+2}{2}}v_{\varepsilon}\\
\leq&(q+\sqrt{n}+1)^2\left\{\eta_2+\eta_2^{-1}q^2(q-2+\sqrt{n})^2\eta_3^\frac{q+2}{q-2}\right\}\int_{\Omega} v_{\varepsilon}^{-q+3}\left|\nabla v_{\varepsilon}\right|^{q-2}\left|D^2 \ln v_{\varepsilon}\right|^2\\
&+\eta_2^{-1}q^2(q-2+\sqrt{n})^2\eta_3^{-\frac{q+2}{4}}\|v_\varepsilon(\cdot,t)\|_{L^\infty(\Omega)}\int_{\Omega}u_\varepsilon^{\frac{q+2}{2}}\\
\end{aligned}
\end{equation}
for all $t>0$ and $\varepsilon\in(0,1)$. Letting $\eta_2=\frac{q(q-1)}{8(q+\sqrt{n}+1)^2}$, $\eta_3=\left(\frac{(q-1)^2}{64(q+\sqrt{n}+1)^4(q-2+\sqrt{n})^2}\right)^{\frac{q-2}{q+2}}$ and using (\ref{Equ(3.9)}), we know that
\begin{equation}\label{Equ(3.10)}
\begin{aligned}
q(q-2+\sqrt{n})\int_{\Omega} u_{\varepsilon} v_{\varepsilon}^{-q+2}&\left|\nabla v_{\varepsilon}\right|^{q-2}\left|D^2 v_{\varepsilon}\right|\leq\frac{q(q-1)}{4}\int_{\Omega} v_{\varepsilon}^{-q+3}\left|\nabla v_{\varepsilon}\right|^{q-2}\left|D^2 \ln v_{\varepsilon}\right|^2\\
&+\frac{8^\frac{q}{2}q(q+\sqrt{n}+1)^q(q-2+\sqrt{n})^\frac{q+2}{2}\|v_\varepsilon(\cdot,t)\|_{L^\infty(\Omega)}}{(q-1)^\frac{q}{2}}
\int_{\Omega}u_\varepsilon^{\frac{q+2}{2}}\\
\end{aligned}
\end{equation}
for all $t>0$ and $\varepsilon\in(0,1)$. Inserting (\ref{Equ(3.8)}) and (\ref{Equ(3.10)}) into (\ref{Equ(3.6)}), we conclude that
\begin{equation}\label{Equ(3.11)}
\begin{aligned}
\frac{d}{dt}&\left(\int_{\Omega}u_\varepsilon^p+\int_{\Omega} v_{\varepsilon}^{-q+1}\left|\nabla v_{\varepsilon}\right|^q\right)+\int_{\Omega}u_\varepsilon^p+\int_{\Omega} v_{\varepsilon}^{-q+1}\left|\nabla v_{\varepsilon}\right|^q\\
&+\frac{q(q-1)}{2}\int_{\Omega} v_{\varepsilon}^{-q+3}\left|\nabla v_{\varepsilon}\right|^{q-2}\left|D^2 \ln v_{\varepsilon}\right|^2\\
\leq& \frac{p(p-1)}{2}\left(\frac{q(q-1)}{2p(p-1)(q+\sqrt{n})^2}\right)^{-\frac{2}{q}}C_2^{\frac{q+2}{q}}
\|v_\varepsilon(\cdot,t)\|_{L^\infty(\Omega)}^\frac{\alpha(q+2)-2}{q}\int_{\Omega}u_\varepsilon^\frac{p(q+2)}{q}\\
&+\frac{8^\frac{q}{2}q(q+\sqrt{n}+1)^q(q-2+\sqrt{n})^\frac{q+2}{2}\|v_\varepsilon(\cdot,t)\|_{L^\infty(\Omega)}}{(q-1)^\frac{q}{2}}
\int_{\Omega}u_\varepsilon^{\frac{q+2}{2}}\\
&+\frac{q}{2} \int_{\partial \Omega} v_{\varepsilon}^{-q+1}\left|\nabla v_{\varepsilon}\right|^{q-2} \cdot \frac{\partial\left|\nabla v_{\varepsilon}\right|^2}{\partial \nu}+\int_{\Omega} v_{\varepsilon}^{-q+1}\left|\nabla v_{\varepsilon}\right|^q+(ap+1)\int_\Omega u_\varepsilon^{p}-bp\int_{\Omega} u_\varepsilon^{p+\gamma-1}
\end{aligned}
\end{equation}
for all $t>0$ and $\varepsilon\in(0,1)$. For any $\eta>0$, an application of (\ref{Equ(2.5)}), Lemmata \ref{result2.3}, \ref{result2.4} and Young's inequality provides some $c_1>0$ such that
\begin{equation}\label{Equ(3.12)}
\begin{aligned}
\frac{q}{2} \int_{\partial \Omega} v_{\varepsilon}^{-q+1}\left|\nabla v_{\varepsilon}\right|^{q-2} \cdot \frac{\partial\left|\nabla v_{\varepsilon}\right|^2}{\partial \nu}\leq&\eta\int_{\Omega}v_{\varepsilon}^{-q-1}\left|\nabla v_{\varepsilon}\right|^{q+2}+\eta\int_{\Omega}v_{\varepsilon}^{-q+1}\left|\nabla v_{\varepsilon}\right|^{q-2}\left|D^2 v_{\varepsilon}\right|^2+c_1\int_{\Omega}v_\varepsilon\\
\leq&2(q+\sqrt{n}+1)^2\eta\int_{\Omega} v_{\varepsilon}^{-q+3}\left|\nabla v_{\varepsilon}\right|^{q-2}\left|D^2 \ln v_{\varepsilon}\right|^2+c_1|\Omega|\|v_0\|_{L^\infty(\Omega)}
\end{aligned}
\end{equation}
and
\begin{equation}\label{Equ(3.13)}
\begin{aligned}
\int_{\Omega} v_{\varepsilon}^{-q+1}\left|\nabla v_{\varepsilon}\right|^q\leq&\eta^\frac{q+2}{q}\int_{\Omega}v_{\varepsilon}^{-q-1}\left|\nabla v_{\varepsilon}\right|^{q+2}+\eta^{-\frac{q+2}{2}}\int_{\Omega}v_\varepsilon\\
\leq&(q+\sqrt{n})^2\eta^\frac{q+2}{q}\int_{\Omega} v_{\varepsilon}^{-q+3}\left|\nabla v_{\varepsilon}\right|^{q-2}\left|D^2 \ln v_{\varepsilon}\right|^2+\eta^{-\frac{q+2}{2}}|\Omega|\|v_0\|_{L^\infty(\Omega)}
\end{aligned}
\end{equation}
for all $t>0$ and $\varepsilon\in(0,1)$. Since $\gamma>1$, by Young's inequality, we can find $c_2>0$ such that
\begin{equation}\label{Equ(3.14)}
\begin{aligned}
(ap+1)\int_\Omega u_\varepsilon^{p}-bp\int_{\Omega} u_\varepsilon^{p+\gamma-1}\leq -\frac{bp}{2}\int_{\Omega} u_\varepsilon^{p+\gamma-1}+c_2\quad \text{for all } t>0 \text{ and  }\varepsilon\in(0,1).
\end{aligned}
\end{equation}
Choosing $\eta$ appropriately small, plugging (\ref{Equ(3.12)})-(\ref{Equ(3.14)}) into (\ref{Equ(3.11)}) we obtain (\ref{Equ(3.5)}). $\hfill{} \Box$

Our next plan is to estimate the first two integral terms on the right-hand sides of (\ref{Equ(3.5)}) appropriately.
\begin{lemma}\label{result3.3}
Let $\alpha\geq1$, $\gamma>1$, $p>1$ and $q\geq2$ be such that
\begin{equation}\label{Equ(3.15)}
\begin{aligned}
q>\frac{2p}{\gamma-1}.
\end{aligned}
\end{equation}
Then there exists a constant $C>0$ such that
\begin{equation}\label{Equ(3.16)}
\begin{aligned}
\frac{p(p-1)}{2}\left(\frac{q(q-1)}{2p(p-1)(q+\sqrt{n})^2}\right)^{-\frac{2}{q}}C_2^{\frac{q+2}{q}}
\|v_0\|_{L^\infty(\Omega)}^\frac{\alpha(q+2)-2}{q}\int_{\Omega}u_\varepsilon^\frac{p(q+2)}{q}
\leq\frac{bp}{4}\int_{\Omega} u_\varepsilon^{p+\gamma-1}+C
\end{aligned}
\end{equation}
for all $t>0$ and $\varepsilon\in(0,1)$, where $C_2$ is given by (\ref{Equ(3.1)}).
\end{lemma}
\noindent{\bf{Proof.}} It follows from (\ref{Equ(3.15)}) that
$$p+\gamma-1-\frac{p(q+2)}{q}=\frac{(\gamma-1)q-2p}{q}>0.$$
Thus, we utilize (\ref{Equ(2.5)}) and Young's inequality to the first summand on the right hand of (\ref{Equ(3.5)}) to show that the existence of $c>0$ such that
\begin{equation*}
\begin{aligned}
\frac{p(p-1)}{2}\left(\frac{q(q-1)}{2p(p-1)(q+\sqrt{n})^2}\right)^{-\frac{2}{q}}C_2^{\frac{q+2}{q}}
\|v_0\|_{L^\infty(\Omega)}^\frac{\alpha(q+2)-2}{q}\int_{\Omega}u_\varepsilon^\frac{p(q+2)}{q}
\leq\frac{bp}{4}\int_{\Omega} u_\varepsilon^{p+\gamma-1}+c
\end{aligned}
\end{equation*}
for all $t>0$ and $\varepsilon\in(0,1)$, which implies (\ref{Equ(3.16)}). $\hfill{} \Box$

\begin{lemma}\label{result3.4}
Let $\gamma>1$, $p>1$ and suppose that $q\geq2$ satisfies
\begin{equation}\label{Equ(3.17)}
\begin{aligned}
q<2(p+\gamma-2).
\end{aligned}
\end{equation}
Then one can find a constant $C>0$ such that
\begin{equation}\label{Equ(3.18)}
\begin{aligned}
\frac{8^\frac{q}{2}q(q+\sqrt{n}+1)^q(q-2+\sqrt{n})^\frac{q+2}{2}\|v_0\|_{L^\infty(\Omega)}}{(q-1)^\frac{q}{2}}
\int_{\Omega}u_\varepsilon^{\frac{q+2}{2}}
\leq\frac{bp}{4}\int_{\Omega} u_\varepsilon^{p+\gamma-1}+C
\end{aligned}
\end{equation}
for all $t>0$ and $\varepsilon\in(0,1)$.
\end{lemma}
\noindent{\bf{Proof.}} Using (\ref{Equ(3.17)}), we obtain
$$p+\gamma-1-\frac{q+2}{2}=p+\gamma-2-\frac{q}{2}>0.$$
An application of Young's inequality and (\ref{Equ(2.5)}) to the second term on the right hand of (\ref{Equ(3.5)}) entails that one can find $c>0$ satisfying
\begin{equation*}
\begin{aligned}
\frac{8^\frac{q}{2}q(q+\sqrt{n}+1)^q(q-2+\sqrt{n})^\frac{q+2}{2}\|v_0\|_{L^\infty(\Omega)}}{(q-1)^\frac{q}{2}}
\int_{\Omega}u_\varepsilon^{\frac{q+2}{2}}
\leq\frac{bp}{4}\int_{\Omega} u_\varepsilon^{p+\gamma-1}+c
\end{aligned}
\end{equation*}
for all $t>0$ and $\varepsilon\in(0,1)$, which immediately gives (\ref{Equ(3.18)}).  $\hfill{} \Box$

In the case $\gamma>2$, as the first application of Lemma \ref{result3.2}, we have the following assertion.
\begin{lemma}\label{result3.5}
Let $\gamma>2$. Assume that (\ref{Equ(1.3)}), (\ref{Equ(1.5)}) and (\ref{Equ(1.6)}) hold with some $\alpha\geq1$. Then for all $p>1$, there exists a constant $C>0$ such that
\begin{equation*}
\|u_\varepsilon(\cdot,t)\|_{L^{p}(\Omega)}\leq C \quad \text{for all } t>0 \text{ and } \varepsilon\in(0,1).
\end{equation*}
\end{lemma}
\noindent{\bf{Proof.}} Since $\gamma>2$, we can readily obtain that
\begin{equation*}
\begin{aligned}
2(p+\gamma-2)-\frac{2p}{\gamma-1}=\frac{2[(p+\gamma-2)(\gamma-1)-p]}{\gamma-1}=\frac{2(\gamma-2)(p+\gamma-1)}{\gamma-1}>0.
\end{aligned}
\end{equation*}
Hence, for any $p >1$, we can pick some $q\geq2$ fulfilling
\begin{equation}\label{Equ(3.19)}
\begin{aligned}
\frac{2p}{\gamma-1}<q<2(p+\gamma-2).
\end{aligned}
\end{equation}
It follows from Lemmata \ref{result3.2}-\ref{result3.4}, (\ref{Equ(2.5)}) and (\ref{Equ(3.19)}) that one can find $c_1>0$ such that
\begin{equation}\label{Equ(3.20)}
\begin{aligned}
\frac{d}{dt}\left(\int_{\Omega}u_\varepsilon^p+\int_{\Omega} v_{\varepsilon}^{-q+1}\left|\nabla v_{\varepsilon}\right|^q\right)+\int_{\Omega}u_\varepsilon^p+\int_{\Omega} v_{\varepsilon}^{-q+1}\left|\nabla v_{\varepsilon}\right|^q
\leq c_1
\end{aligned}
\end{equation}
for all $t>0$ and $\varepsilon\in(0,1)$. Letting $y_\varepsilon(t):=\int_{\Omega}u_\varepsilon^p(\cdot,t)+\int_{\Omega}v^{-q+1}_{\varepsilon}(\cdot,t)\left|\nabla v_{\varepsilon}(\cdot,t)\right|^q$ and using (\ref{Equ(3.20)}), we obtain $y'_\varepsilon(t)+y_\varepsilon(t)\leq c_1$ for all $t>0$ and $\varepsilon\in(0,1)$, thus an ODE argument implies the boundedness of $\|u_\varepsilon(\cdot,t)\|_{L^{p}(\Omega)}$  for all $t>0$ and $\varepsilon\in(0,1)$. $\hfill{} \Box$

In the case $\gamma=2$ and $n\geq3$, we will derive $L^p$-bounds on $u_\varepsilon$ if $b$ is sufficiently large.
\begin{lemma}\label{result3.6}
Let $\gamma=2$, $n\geq3$ and $p>1$. Suppose that (\ref{Equ(1.3)}), (\ref{Equ(1.5)}) and (\ref{Equ(1.6)}) hold, and that $b>0$ satisfies
\begin{equation}\label{Equ(3.21)}
\begin{aligned}
b>\kappa_1(p,n)C_2^{\frac{p+1}{p}}\|v_0\|_{L^\infty(\Omega)}^{\frac{\alpha(p+1)-1}{p}}+\kappa_2(p,n)\|v_0\|_{L^\infty(\Omega)},
\end{aligned}
\end{equation}
where $C_2$ is as defined in Lemma \ref{result3.1},
\begin{equation}\label{Equ(3.22)}
\begin{split}
\kappa_1(p,n):=(p-1)^{\frac{p+1}{p}}(2p+\sqrt{n})^\frac{2}{p}(2p-1)^{-\frac{1}{p}}
\end{split}
\end{equation}
and
\begin{equation}\label{Equ(3.23)}
\begin{split}
\kappa_2(p,n):=\frac{2^{3p+2}(2p+\sqrt{n}+1)^{2p}(2p-2+\sqrt{n})^{p+1}}{(2p-1)^p}.
\end{split}
\end{equation}
Then we infer the existence of $C>0$ such that
\begin{equation}\label{Equ(3.24)}
\|u_\varepsilon(\cdot,t)\|_{L^p(\Omega)}\leq C \quad \text{for all } t>0 \text{ and } \varepsilon\in(0,1).
\end{equation}
\end{lemma}
\noindent{\bf{Proof.}} Since $\gamma=2$, using (\ref{Equ(2.5)}), (\ref{Equ(3.22)}) and (\ref{Equ(3.23)}), we apply Lemma \ref{result3.2} to $q:=2p$ to find $c>0$ satisfying
\begin{equation}\label{Equ(3.25)}
\begin{aligned}
\frac{d}{dt}&\left(\int_{\Omega}u_\varepsilon^p+\int_{\Omega} v_{\varepsilon}^{-2p+1}\left|\nabla v_{\varepsilon}\right|^{2p}\right)+\int_{\Omega}u_\varepsilon^p+\int_{\Omega} v_{\varepsilon}^{-2p+1}\left|\nabla v_{\varepsilon}\right|^{2p}\\
\leq& \frac{p(p-1)}{2}\left(\frac{2p-1}{(p-1)(2p+\sqrt{n})^2}\right)^{-\frac{1}{p}}C_2^{\frac{p+1}{p}}
\|v_0\|_{L^\infty(\Omega)}^\frac{\alpha(p+1)-1}{p}\int_{\Omega}u_\varepsilon^{p+1}\\
&+\frac{2^{3p+1}p(2p+\sqrt{n}+1)^{2p}(2p-2+\sqrt{n})^{p+1}\|v_0\|_{L^\infty(\Omega)}}{(2p-1)^p}
\int_{\Omega}u_\varepsilon^{p+1}-\frac{bp}{2}\int_{\Omega} u_\varepsilon^{p+1}+c\\
=& -\frac{p}{2}\left\{b-\kappa_1(p,n)C_2^{\frac{p+1}{p}}
\|v_0\|_{L^\infty(\Omega)}^\frac{\alpha(p+1)-1}{p}-\kappa_2(p,n)\|v_0\|_{L^\infty(\Omega)}\right\}\int_{\Omega}u_\varepsilon^{p+1}+c
\end{aligned}
\end{equation}
for all $t>0$ and $\varepsilon\in(0,1)$. Letting $y_\varepsilon(t):=\int_{\Omega}u_\varepsilon^p(\cdot,t)+\int_{\Omega}v^{-2p+1}_{\varepsilon}(\cdot,t)\left|\nabla v_{\varepsilon}(\cdot,t)\right|^{2p}$, applying (\ref{Equ(3.25)}) and (\ref{Equ(3.21)}), we obtain $y'_\varepsilon(t)+y_\varepsilon(t)\leq c$ for all $t>0$ and $\varepsilon\in(0,1)$. Therefore, an ODE argument shows (\ref{Equ(3.24)}). $\hfill{} \Box$

While for the lower-dimensional settings in the case $\gamma=2$, inspired by an idea in \cite{WAA}, the $L^2$ estimates of $u_\varepsilon$ can be established for arbitrary $b>0$.
\begin{lemma}\label{result3.7}
Let $n\leq2$ and $\gamma=2$, and suppose that (\ref{Equ(1.3)}), (\ref{Equ(1.5)}) and (\ref{Equ(1.6)}) with $\alpha\geq1$ hold. Then there exists a constant $C>0$ such that
\begin{equation*}
\|u_\varepsilon(\cdot,t)\|_{L^2(\Omega)}\leq C \quad \text{for all } t>0 \text{ and } \varepsilon\in(0,1).
\end{equation*}
\end{lemma}
\noindent{\bf{Proof.}} It follows from $n\leq2$ that $W^{1,1}(\Omega)$ is continuously embedded into $L^2(\Omega)$. We thus utilize a Sobolev inequality, H\"{o}lder's and Young's inequalities to show that
\begin{equation}\label{Equ(3.26)}
\|\psi\|_{L^2(\Omega)}\leq c_1\|\nabla\psi\|_{L^1(\Omega)}+c_1\|\psi\|_{L^\frac{1}{2}(\Omega)} \quad \text{for all } \psi\in W^{1,1}(\Omega)
\end{equation}
with some $c_1>0$ (see \cite[Lemma 3.5]{WAA} for details). Lemma \ref{result3.1} can be applied to $p:=2$ to obtain that one can find some $c_2>0$ and $c_3>0$ such that
\begin{equation}\label{Equ(3.27)}
\begin{aligned}
&\frac{d}{dt}\int_{\Omega} u_\varepsilon^{2}+c_2\int_\Omega v_\varepsilon^\alpha|\nabla u_\varepsilon|^2+2b\int_{\Omega} u_\varepsilon^3\leq c_3\int_{\Omega} u_\varepsilon^2v_\varepsilon^{\alpha-2}|\nabla v_\varepsilon|^2+
2a\int_\Omega u_\varepsilon^{2}
\end{aligned}
\end{equation}
for all $t>0$ and $\varepsilon\in(0,1)$. Due to $\alpha\geq1$, we use (\ref{Equ(2.5)}) and H\"{o}lder's inequality to conclude
\begin{equation}\label{Equ(3.28)}
\begin{aligned}
\int_{\Omega} u_\varepsilon^2v_\varepsilon^{\alpha-2}|\nabla v_\varepsilon|^2\leq \left(\int_\Omega \frac{|\nabla v_\varepsilon|^4}{v^3_\varepsilon}\right)^\frac{1}{2}\left(\int_\Omega u_\varepsilon^4v^{2\alpha-1}_\varepsilon\right)^\frac{1}{2}\leq \|v_0\|^\frac{\alpha-1}{2}_{L^{\infty}(\Omega)}\left(\int_\Omega \frac{|\nabla v_\varepsilon|^4}{v^3_\varepsilon}\right)^\frac{1}{2}\left(\int_\Omega u_\varepsilon^4v^\alpha_\varepsilon\right)^\frac{1}{2}
\end{aligned}
\end{equation}
for all $t>0$ and $\varepsilon\in(0,1)$. An application of  (\ref{Equ(3.26)}), (\ref{Equ(2.3)}), (\ref{Equ(2.5)}) and H\"{o}lder's inequality implies
\begin{equation}\label{Equ(3.29)}
\begin{aligned}
\left(\int_\Omega u_\varepsilon^4v^\alpha_\varepsilon\right)^\frac{1}{2}=&\|u_\varepsilon^2v^\frac{\alpha}{2}_\varepsilon\|_{L^2(\Omega)}\\
\leq&2c_{1}\int_\Omega u_\varepsilon v^\frac{\alpha}{2}_\varepsilon|\nabla u_\varepsilon|+\frac{c_{1}\alpha}{2}\int_\Omega u_\varepsilon^2v^\frac{\alpha-2}{2}_\varepsilon|\nabla v_\varepsilon|+c_1\left(\int_\Omega u_\varepsilon v^\frac{\alpha}{4}_\varepsilon\right)^2\\
\leq&2c_{1}\left(\int_\Omega u_\varepsilon^2\right)^\frac{1}{2}\left(\int_\Omega v_\varepsilon^\alpha|\nabla u_\varepsilon|^2\right)^\frac{1}{2}+\frac{c_{1}\alpha}{2}\left(\int_\Omega u_\varepsilon^2\right)^\frac{1}{2}\left(\int_\Omega u_\varepsilon^2v_\varepsilon^{\alpha-2}|\nabla v_\varepsilon|^2\right)^\frac{1}{2}\\
&+c_1m_1^2\|v_0\|^\frac{\alpha}{2}_{L^{\infty}(\Omega)}\\
\end{aligned}
\end{equation}
for all $t>0$ and $\varepsilon\in(0,1)$ with $m_1$ taken from (\ref{Equ(2.3)}). We employ (\ref{Equ(3.28)}), (\ref{Equ(3.29)}) and Young's inequality to show that
\begin{equation*}
\begin{aligned}
c_3\int_{\Omega} u_\varepsilon^2v_\varepsilon^{\alpha-2}|\nabla v_\varepsilon|^2\leq&c_3\|v_0\|^\frac{\alpha-1}{2}_{L^{\infty}(\Omega)}\left(\int_\Omega \frac{|\nabla v_\varepsilon|^4}{v^3_\varepsilon}\right)^\frac{1}{2}\Bigg[2c_{1}\left(\int_\Omega u_\varepsilon^2\right)^\frac{1}{2}\left(\int_\Omega v_\varepsilon^\alpha|\nabla u_\varepsilon|^2\right)^\frac{1}{2}\\
&+\frac{c_{1}\alpha}{2}\left(\int_\Omega u_\varepsilon^2\right)^\frac{1}{2}\left(\int_\Omega u_\varepsilon^2v_\varepsilon^{\alpha-2}|\nabla v_\varepsilon|^2\right)^\frac{1}{2}
+c_1m_1^2\|v_0\|^\frac{\alpha}{2}_{L^{\infty}(\Omega)}\Bigg]\\
\leq&\frac{c_2}{2}\int_\Omega v_\varepsilon^\alpha|\nabla u_\varepsilon|^2+\frac{2c_1^2c^2_3\|v_0\|^{\alpha-1}_{L^{\infty}(\Omega)}}{c_2}\left(\int_\Omega u_\varepsilon^2\right)\left(\int_\Omega \frac{|\nabla v_\varepsilon|^4}{v^3_\varepsilon}\right)\\
&+\frac{c_3}{2}\int_{\Omega} u_\varepsilon^2v_\varepsilon^{\alpha-2}|\nabla v_\varepsilon|^2+\frac{c_1^2c_3\alpha^2\|v_0\|^{\alpha-1}_{L^{\infty}(\Omega)}}{8}\left(\int_\Omega u_\varepsilon^2\right)\left(\int_\Omega \frac{|\nabla v_\varepsilon|^4}{v^3_\varepsilon}\right)\\
&+\frac{1}{2}\int_\Omega \frac{|\nabla v_\varepsilon|^4}{v^3_\varepsilon}+\frac{c_1^2c_3^2m_1^4\|v_0\|^{2\alpha-1}_{L^{\infty}(\Omega)}}{2}
\end{aligned}
\end{equation*}
for all $t>0$ and $\varepsilon\in(0,1)$, which immediately implies
\begin{equation}\label{Equ(3.30)}
\begin{aligned}
c_3\int_{\Omega} u_\varepsilon^2v_\varepsilon^{\alpha-2}|\nabla v_\varepsilon|^2
\leq&c_2\int_\Omega v_\varepsilon^\alpha|\nabla u_\varepsilon|^2+\int_\Omega \frac{|\nabla v_\varepsilon|^4}{v^3_\varepsilon}+c_1^2c_3^2m_1^4\|v_0\|^{2\alpha-1}_{L^{\infty}(\Omega)}\\
&+\frac{c_1^2c_3\|v_0\|^{\alpha-1}_{L^{\infty}(\Omega)}}{4}\left(\alpha^2+\frac{16c_3}{c_2}\right)\left(\int_\Omega u_\varepsilon^2\right)\left(\int_\Omega \frac{|\nabla v_\varepsilon|^4}{v^3_\varepsilon}\right)
\end{aligned}
\end{equation}
for all $t>0$ and $\varepsilon\in(0,1)$. Using (\ref{Equ(3.27)}), (\ref{Equ(3.30)}) and H\"{o}lder's inequality, we see that
\begin{equation}\label{Equ(3.31)}
\begin{aligned}
\frac{d}{dt}\int_{\Omega} u_\varepsilon^{2}+2b\int_{\Omega} u_\varepsilon^3
\leq&\frac{c_1^2c_3\|v_0\|^{\alpha-1}_{L^{\infty}(\Omega)}}{4}\left(\alpha^2+\frac{16c_3}{c_2}\right)\left(\int_\Omega u_\varepsilon^2\right)\left(\int_\Omega \frac{|\nabla v_\varepsilon|^4}{v^3_\varepsilon}\right)\\
&+\int_\Omega \frac{|\nabla v_\varepsilon|^4}{v^3_\varepsilon}+2a\int_\Omega u_\varepsilon^{2}+c_1^2c_3^2m_1^4\|v_0\|^{2\alpha-1}_{L^{\infty}(\Omega)}
\end{aligned}
\end{equation}
for all $t>0$ and $\varepsilon\in(0,1)$. By  Young's and H\"{o}lder's inequalities we know that
\begin{equation}\label{Equ(3.32)}
\begin{aligned}
2^\frac{1}{2}b|\Omega|^{-\frac{1}{2}}\left(\int_{\Omega} u_\varepsilon^2+1\right)^\frac{3}{2}\leq2b|\Omega|^{-\frac{1}{2}}\left(\int_{\Omega} u_\varepsilon^2\right)^\frac{3}{2}+2b|\Omega|^{-\frac{1}{2}}\leq2b\int_{\Omega} u_\varepsilon^3+2b|\Omega|^{-\frac{1}{2}}
\end{aligned}
\end{equation}
for all $t>0$ and $\varepsilon\in(0,1)$. Hence, it follows from (\ref{Equ(3.31)}) and (\ref{Equ(3.32)}) that one can find $c_4>0$ and $c_5>0$ fulfilling
\begin{equation}\label{Equ(3.33)}
\begin{aligned}
\frac{d}{dt}\int_{\Omega} u_\varepsilon^{2}+2^\frac{1}{2}b|\Omega|^{-\frac{1}{2}}\left(\int_{\Omega} u_\varepsilon^2+1\right)^\frac{3}{2}\leq\left(\int_{\Omega} u_\varepsilon^{2}+1\right)\left(c_4\int_{\Omega} \frac{|\nabla v_\varepsilon|^4}{v^3_\varepsilon}+c_5\right)
\end{aligned}
\end{equation}
for all $t>0$ and $\varepsilon\in(0,1)$. Now, for each $\varepsilon\in(0,1)$, letting $y_\varepsilon(t):=\int_{\Omega} u_\varepsilon^{2}(\cdot,t)+1$ for all $t\geq0$ and $g_\varepsilon(t):=c_4\int_{\Omega} \frac{|\nabla v_\varepsilon(\cdot,t)|^4}{v^3_\varepsilon(\cdot,t)}+c_5$ for all $t>0$, (\ref{Equ(3.33)}) shows
$$y_\varepsilon'(t)+2^\frac{1}{2}b|\Omega|^{-\frac{1}{2}}y^\frac{3}{2}_\varepsilon(t)\leq y_\varepsilon(t)g_\varepsilon(t) \quad \text{for all } t>0 \text{ and } \varepsilon\in(0,1).$$
This along with (\ref{Equ(2.15)}) and Lemma \ref{result2.7} entails the intended conclusion. $\hfill{} \Box$

Although we achieve the $L^p$-bounds for $u_\varepsilon$, unfortunately, the uniform boundedness of solutions can not be obtained through the outcome of a Moser-type iterative, because of the possible degeneracy of diffusion for the first equation in (\ref{Equ(2.1)}). However, by means of a transformation $w_{\varepsilon}=-\ln \frac{v_{\varepsilon}}{\left\|v_0\right\|_{L^{\infty}(\Omega)}}$ dating back to \cite{WM3}, we can construct a time-dependent pointwise lower bound for $v_\varepsilon$ to obtain local-in-time $L^\infty$-estimates for $u_\varepsilon$.
\begin{lemma}\label{result3.8}
Let (\ref{Equ(1.3)}) hold and $n\geq1$. For any $p>\frac{n}{2}$, assume that there exists $C_1(p)>0$  such that
\begin{equation}\label{Equ(3.34)}
\begin{aligned}
\|u_\varepsilon(\cdot,t)\|_{L^p(\Omega)}\leq C_1(p) \quad \text{for all } t>0 \text{ and } \varepsilon\in(0,1).
\end{aligned}
\end{equation}
Then given any $T>0$ there exists $C_2(T,p)>0$ such that
\begin{equation}\label{Equ(3.35)}
v_\varepsilon(x,t)\geq C_2(T,p) \quad \text{for all } x\in\Omega, \,\,t\in(0,T) \text{ and } \varepsilon\in(0,1).
\end{equation}
\end{lemma}
\noindent{\bf{Proof.}} Let $w_{\varepsilon}(x, t):=-\ln \frac{v_{\varepsilon}(x, t)}{\left\|v_0\right\|_{L^{\infty}(\Omega)}}$ in $(x,t)\in\overline{\Omega}\times[0,\infty)$ for all $\varepsilon\in(0,1).$ Then for each $\varepsilon\in(0,1)$, we see that $w_{\varepsilon}\geq0$ in $\Omega\times(0,\infty)$ and that
\begin{equation}\label{Equ(3.36)}
\begin{cases}
w_{\varepsilon t}=\Delta w_{\varepsilon}-|\nabla w_{\varepsilon}|^2+\frac{u_{\varepsilon}}{1+\varepsilon u_{\varepsilon}}, & x \in \Omega,\,\,t>0, \\
\frac{\partial w_{\varepsilon}}{\partial \nu}=0, & x \in \partial \Omega,\,\,t>0, \\
w_{\varepsilon}(x, 0)=w_0(x):=-\ln \frac{v_0(x)}{\left\|v_0\right\|_{L^{\infty}(\Omega)}}, & x \in \Omega.
\end{cases}
\end{equation}
On the basis of the variation-of-constants formula of (\ref{Equ(3.36)}), one has
\begin{equation*}
\begin{aligned}
w_\varepsilon(\cdot,t)=&e^{t\Delta}w_0-\int_{0}^{t}e^{(t-s)\Delta}|\nabla w_{\varepsilon}(\cdot,s)|^2+\int_{0}^{t}e^{(t-s)\Delta}\frac{u_\varepsilon(\cdot,s)}{1+\varepsilon u_\varepsilon(\cdot,s)}\\
\leq&e^{t\Delta}w_0+\int_{0}^{t}e^{(t-s)\Delta}\frac{u_\varepsilon(\cdot,s)}{1+\varepsilon u_\varepsilon(\cdot,s)}\\
\end{aligned}
\end{equation*}
for all $t>0$ and $\varepsilon\in (0,1)$. Then by virtue of (\ref{Equ(3.34)}) and the smoothing properties of Neumann heat semigroup $\left(e^{t \Delta}\right)_{t \geq 0}$ on $\Omega$ (\cite[Lemma 1.3]{WJ}), we infer the existences of $c_1>0$ and $c_2>0$ such that
\begin{equation}\label{Equ(3.37)}
\begin{aligned}
\|w_\varepsilon(\cdot,t)\|_{L^{\infty}(\Omega)}\leq& \|e^{t\Delta}w_0\|_{L^{\infty}(\Omega)}+\int_{0}^{t} \| e^{(t-s)\Delta}\frac{u_\varepsilon(\cdot,s)}{1+\varepsilon u_\varepsilon(\cdot,s)}\|_{L^{\infty}(\Omega)}ds\\
\leq& \|w_0\|_{L^{\infty}(\Omega)}+c_1\int_{0}^{t} \{1+(t-s)^{-\frac{n}{2p}}\}\|u_\varepsilon(\cdot,s)\|_{L^{p}(\Omega)}ds\\
\leq& \|w_0\|_{L^{\infty}(\Omega)}+c_2\int_{0}^{t} (1+\sigma^{-\frac{n}{2p}})d\sigma\\
\end{aligned}
\end{equation}
for all $t>0$ and $\varepsilon\in (0,1)$. It follows from $p>\frac{n}{2}$ that
\begin{equation}\label{Equ(3.38)}
\begin{aligned}
\int_{0}^{t} (1+\sigma^{-\frac{n}{2p}})d\sigma=t+\frac{2p}{2p-n}t^\frac{2p-n}{2p} \text{ for all } t>0.
\end{aligned}
\end{equation}
Hence, for any $T>0$, (\ref{Equ(3.37)}) in conjunction with (\ref{Equ(3.38)}) implies that one can find some $c_3(T,p)>0$ such that
$$\|w_\varepsilon(\cdot,t)\|_{L^{\infty}(\Omega)}\leq c_3(T,p)\quad \text{for all } t\in(0,T) \text{ and } \varepsilon\in(0,1).$$
Then according to the definition of $w_{\varepsilon}(x, t)$, we can readily obtain (\ref{Equ(3.35)}). $\hfill{} \Box$

With the lower bound of $v_\varepsilon$ at hand, we can show local boundedness criterion of solutions of (\ref{Equ(2.1)}).
\begin{lemma}\label{result3.9}
Let (\ref{Equ(1.3)}), (\ref{Equ(1.5)}) and (\ref{Equ(1.6)}) hold. For all $n\geq1$, assume that there exist $C>0$ and $p\geq1$ such that $p>\frac{n}{2}$ and
\begin{equation}\label{Equ(3.39)}
\begin{aligned}
\|u_\varepsilon(\cdot,t)\|_{L^p(\Omega)}\leq C \quad \text{for all } t>0 \text{ and } \varepsilon\in(0,1).
\end{aligned}
\end{equation}
Then for all $T>0$, one can find $C(T)>0$ such that
\begin{equation*}
\|u_\varepsilon(\cdot,t)\|_{L^\infty(\Omega)}\leq C(T) \quad \text{for all } t\in(0,T)\text{ and } \varepsilon\in(0,1).
\end{equation*}
\end{lemma}
\noindent{\bf{Proof.}}  Without loss of generality, we suppose that $p\leq n$. Using $\frac{n}{2} < p\leq n$, we readily obtain $\frac{np}{n - p}> n$. Hence, together with $p\geq1$, it follows that we can find $r>2$ such that $n<r<\frac{np}{n-p}$. An application of (\ref{Equ(3.39)}) and a regularization features of the Neumann heat semigroup (see \cite[Lemma 4.1]{DW}) implies that there exists $c_1>0$ such that $\|\nabla v_\varepsilon(\cdot,t)\|_{L^{r}(\Omega)}\leq c_1$ for all $t>0$ and $\varepsilon\in(0,1)$. Moreover, (\ref{Equ(3.39)}) in conjunction with Lemma \ref{result3.8} shows that for any $T>0$, we can find $c_2(T)>0$ fulfilling $v_\varepsilon\geq c_2(T)$ in $\Omega\times(0,T)$ for all $\varepsilon\in(0,1)$. Similar to (\ref{Equ(3.1)}), there exist $c_3>0$ and $c_4>0$ such that for all $q>1$ we have
\begin{equation}\label{Equ(3.40)}
\begin{aligned}
\frac{d}{dt}\int_{\Omega} u_\varepsilon^{q}+\frac{2 c_3c^\alpha_2(T)(q-1)}{q}\int_\Omega|\nabla u^\frac{q}{2}_\varepsilon|^2+bq\int_{\Omega} u_\varepsilon^{q+\gamma-1}\leq& \frac{q(q-1)c_4 c_5(T)}{2}\int_{\Omega} u_\varepsilon^{q}|\nabla v_\varepsilon|^2\\
&+aq\int_\Omega u_\varepsilon^{q}
\end{aligned}
\end{equation}
for all $ t\in(0,T)$ and $\varepsilon\in(0,1)$, where $c_5(T):=\max\{c^{\alpha-2}_2(T),\|v_0\|^{\alpha-2}_{L^\infty(\Omega)}\}$. An application of $r>\max\{2,n\}$ entails $\frac{2r}{r-2}<\frac{2n}{(n-2)_+}$, thus $W^{1,2}(\Omega)$ is compactly embedded into $L^\frac{2r}{r-2}(\Omega)$ and into $L^2(\Omega)$. Applying H\"{o}lder's inequality and an Ehrling type inequality and (\ref{Equ(2.3)}), we infer the existence of $c_6(q,T)>0$ satisfying
\begin{equation}\label{Equ(3.41)}
\begin{aligned}
\frac{q(q-1)c_4 c_5(T)}{2}\int_{\Omega} u_\varepsilon^{q}|\nabla v_\varepsilon|^2+aq\int_\Omega u_\varepsilon^{q}
\leq&\frac{q(q-1)c_4 c_5(T)}{2}\|u_\varepsilon^\frac{q}{2}\|^2_{L^\frac{2r}{r-2}(\Omega)}\|\nabla v_\varepsilon\|^2_{L^{r}(\Omega)}\\
&+aq\|u_\varepsilon^\frac{q}{2}\|^2_{L^2(\Omega)}\\
\leq&\frac{q(q-1)c_1^2c_4c_5(T)}{2}\|u_\varepsilon^\frac{q}{2}\|^2_{L^\frac{2r}{r-2}(\Omega)}+aq\|u_\varepsilon^\frac{q}{2}\|^2_{L^2(\Omega)}\\
\leq&\frac{2 c_3 c^\alpha_2(T)(q-1)}{q}\int_\Omega |\nabla u_\varepsilon^\frac{q}{2}|^2+c_6(q,T)
\end{aligned}
\end{equation}
for all $ t \in (0,T)$ and $\varepsilon\in(0,1)$. We utilize $\gamma>1$ and H\"{o}lder's inequality to obtain
\begin{equation}\label{Equ(3.42)}
\begin{aligned}
bq|\Omega|^\frac{1-\gamma}{q}\left(\int_{\Omega} u_\varepsilon^{q}\right)^\frac{q+\gamma-1}{q}\leq bq\int_{\Omega} u_\varepsilon^{q+\gamma-1}\quad \text{for all } t>0 \text{ and } \varepsilon\in(0,1).
\end{aligned}
\end{equation}
Letting $y_\varepsilon(t):=\int_{\Omega}u_\varepsilon^q(\cdot,t)$ for all $ t>0$, $\varepsilon\in(0,1)$ and collecting (\ref{Equ(3.40)})-(\ref{Equ(3.42)}), we see that $$y_\varepsilon'(t)+bq|\Omega|^\frac{1-\gamma}{q}y_\varepsilon^\frac{q+\gamma-1}{q}(t)\leq c_6(q,T) \text{ for all } t\in (0,T) \text{ and } \varepsilon\in(0,1).$$
We thus employ an ODE argument to find some $c_7(q,T)>0$ such that $\|u_\varepsilon(\cdot,t)\|_{L^q(\Omega)}\leq c_7(q,T)$ for all $t\in (0,T)$ and $\varepsilon\in(0,1)$. Thus, the regularization features of the Neumann heat semigroup entail that for some $c_8(q,T)>0$, we have $\|v_\varepsilon(\cdot,t)\|_{W^{1,\infty}(\Omega)}\leq c_8(q,T)$ for all $t\in (0,T)$ and $\varepsilon\in(0,1)$. Then by Moser-iteration procedure \cite[Lemma A.1 ]{TWJ} we deduce the existence of $c_9(T)>0$ such that $\|u_\varepsilon(\cdot,t)\|_{L^\infty(\Omega)}\leq c_9(T)$ for all $t\in (0,T)$ and $\varepsilon\in(0,1)$. Hence, we can derive the claimed conclusion. $\hfill{} \Box$

With the aid of Lemma \ref{result3.9}, the following result can assert that the global solutions are locally bounded under some assumptions.
\begin{lemma}\label{result3.10}
Suppose that (\ref{Equ(1.3)}), (\ref{Equ(1.5)}) and (\ref{Equ(1.6)}) with $\alpha\geq1$ hold. For all $n\geq1$, if one of the following cases holds:\\
(i) $\gamma>2$;\\
(ii) $\gamma=2$ and $b>\left\{
\begin{array}{llll}
0,\quad &\text{if } n\leq2,\\
\kappa_1\left(\frac{n}{2},n\right)C_2^{\frac{n+2}{n}}\|v_0\|_{L^\infty(\Omega)}^{\frac{(n+2)\alpha-2}{n}}+\kappa_2\left(\frac{n}{2},n\right)\|v_0\|_{L^\infty(\Omega)},\ &\text{if }n\geq3,
\end{array}
\right.$\\
where $\kappa_1$ and $\kappa_2$ are as defined in (\ref{Equ(3.22)}) and (\ref{Equ(3.23)}), respectively, and $C_2$ is given by Lemma \ref{result3.1}, then for all $T>0$, there exists $C(T)>0$ such that
\begin{equation*}
\|u_\varepsilon(\cdot,t)\|_{L^\infty(\Omega)}+\|v_\varepsilon(\cdot,t)\|_{W^{1,\infty}(\Omega)}\leq C(T) \quad \text{for all } t\in(0,T)\text{ and } \varepsilon\in(0,1).
\end{equation*}
\end{lemma}
\noindent{\bf{Proof.}} Let $p\geq1$ such that $p>\frac{n}{2}$. In the case $\gamma>2$, Lemma \ref{result3.5} guarantees that $\|u_\varepsilon(\cdot,t)\|_{L^{p}(\Omega)}$ is bounded for all $t>0$ and $\varepsilon\in(0,1)$. In the case $\gamma=2$ and $n\leq2$, for any $b>0$, the uniform boundedness of $\|u_\varepsilon(\cdot,t)\|_{L^2(\Omega)}$ for all $t>0$ and $\varepsilon \in(0,1)$ is proved in Lemma \ref{result3.7}. Whereas in the case $\gamma=2$ and $n\geq3$, thanks to
\begin{equation*}
\begin{aligned}
b>\kappa_1\left(\frac{n}{2},n\right)C_2^{\frac{n+2}{n}}\|v_0\|_{L^\infty(\Omega)}^
{\frac{\alpha(n+2)-2}{n}}+\kappa_2\left(\frac{n}{2},n\right)\|v_0\|_{L^\infty(\Omega)}
\end{aligned}
\end{equation*}
and the continuity of $\kappa_1$ and $\kappa_2$, there exists $p>\frac{n}{2}$  satisfying
\begin{equation*}
\begin{aligned}
b>\kappa_1(p,n)C_2^{\frac{p+1}{p}}\|v_0\|_{L^\infty(\Omega)}^{\frac{\alpha(p+1)-1}{p}}+\kappa_2(p,n)\|v_0\|_{L^\infty(\Omega)}.
\end{aligned}
\end{equation*}
Applying Lemma \ref{result3.6} we also see that $\|u_\varepsilon(\cdot,t)\|_{L^{p}(\Omega)}$ is bounded for all $t>0$ and $\varepsilon\in(0,1)$. In conclusion, we utilize Lemma \ref{result3.9} to complete this proof. $\hfill{} \Box$

As a straightforward consequence of Lemma \ref{result3.10}, we can obtain H\"{o}der regularity properties of $u_\varepsilon$ and $v_\varepsilon$.
\begin{lemma}\label{result3.11}
Let the assumptions of Theorem \ref{result1.1} hold. For any $T>0$, there exist  $C(T)>0$ and $\sigma=\sigma(T)\in(0,1)$ such that
\begin{equation*}
\|u_\varepsilon\|_{C^{\sigma, \frac{\sigma}{2}}(\bar{\Omega} \times[0, T])}+\|v_\varepsilon\|_{C^{\sigma, \frac{\sigma}{2}}(\bar{\Omega} \times[0, T])} \leq C(T) \quad \text{for all } \varepsilon \in(0,1).
\end{equation*}
\end{lemma}
\noindent{\bf{Proof.}} We rewrite the first equation of model (\ref{Equ(2.1)}) as follows
\begin{equation*}
u_{\varepsilon t}=\nabla\cdot A_\varepsilon(x,t,u_\varepsilon,\nabla u_\varepsilon)+B_\varepsilon(x,t,u_\varepsilon)\quad \,\,x\in\Omega,\,\,t>0,\,\, \varepsilon \in(0,1)
\end{equation*}
with $$A_\varepsilon(x,t,u_\varepsilon,\nabla u_\varepsilon):=\varepsilon \nabla u_\varepsilon(x,t)+\phi(v_\varepsilon(x,t))\nabla u_\varepsilon(x,t)+\phi'(v_\varepsilon(x,t))u_\varepsilon(x,t)\nabla v_\varepsilon(x,t)$$
and
$$B_\varepsilon(x,t,u_\varepsilon):=au_\varepsilon(x,t)-bu^\gamma_\varepsilon(x,t).$$
An application of Lemmata \ref{result3.10}, \ref{result3.8}, (\ref{Equ(1.4)}) and (\ref{Equ(2.5)}) entails that there exist some $c_i(T)>0$ ($i=1,2,3,4$) such that
$$A_\varepsilon(x,t,u_\varepsilon,\nabla u_\varepsilon)\nabla u_\varepsilon\geq c_1(T)|\nabla u_\varepsilon|^2-c_2(T)\quad \,\,x\in\Omega,\,\,t\in(0,T),\,\, \varepsilon \in(0,1)$$
and
$$|A_\varepsilon(x,t,u_\varepsilon,\nabla u_\varepsilon)|\leq c_3(T)|\nabla u_\varepsilon|+c_3(T)\quad \,\,x\in\Omega,\,\,t\in(0,T),\,\, \varepsilon \in(0,1)$$
as well as
$$|B_\varepsilon(x,t,u_\varepsilon)|\leq c_4(T)\quad \,\,x\in\Omega,\,\,t\in(0,T),\,\, \varepsilon \in(0,1).$$
Therefore, this in conjunction with (\ref{Equ(1.3)}) and \cite[Theorem 1.3 and Remark 1.4]{PV} implies that one can find some $c_5(T)>0$ and $\sigma_1=\sigma_1(T)\in(0,1)$ such that
$\|u_\varepsilon\|_{C^{\sigma_1, \frac{\sigma_1}{2}}(\bar{\Omega} \times[0, T])} \leq c_5(T)$ for all  $\varepsilon \in(0,1).$
Similarly, from the second equation of model (\ref{Equ(2.1)}), we have $\|v_\varepsilon\|_{C^{\sigma_2, \frac{\sigma_2}{2}}(\bar{\Omega} \times[0, T])} \leq c_6(T)$ with some $c_6(T)>0$ and $\sigma_2=\sigma_2(T)\in(0,1)$ for all  $\varepsilon \in(0,1).$ Hence the proof of Lemma \ref{result3.11} is completed. $\hfill{} \Box$

An application of Lemma \ref{result3.11} and parabolic Schauder theory leads to the following higher order estimates.
\begin{lemma}\label{result3.12}
Under the assumptions of Theorem \ref{result1.1}, for all $T>0$ and any $\tau\in(0,T)$, there exist  $C(T,\tau)>0$ and $\sigma=\sigma(T,\tau)\in(0,1)$ such that
\begin{equation}\label{Equ(3.43)}
\|u_\varepsilon\|_{C^{2+\sigma, 1+\frac{\sigma}{2}}(\bar{\Omega} \times[\tau, T])}+\|v_\varepsilon\|_{C^{2+\sigma, 1+\frac{\sigma}{2}}(\bar{\Omega} \times[\tau, T])} \leq C(T,\tau) \,\,\text{for all } \varepsilon \in(0,1).
\end{equation}
\end{lemma}
\noindent{\bf{Proof.}} Along with Lemma \ref{result3.11} and \cite[Theorem IV.5.3]{LSU} we infer that (\ref{Equ(3.43)}) is valid, for more details, we refer to \cite[Lemma 4.5 ]{WAA}.  $\hfill{} \Box$

Now, we begin with the proof of Theorem \ref{result1.1}. \\
{\bf Proof of Theorem \ref{result1.1}.}  As a consequence of Lemmata \ref{result3.11}, \ref{result3.12} and the Arzel\`{a}-Ascoli theorem,
there exist a sequence $\left(\varepsilon_{j}\right)_{j \in \mathbb{N}} \subset(0,1)$ and functions $u,v$ such that $\varepsilon_{j} \searrow 0$ as $j \rightarrow \infty$ and that
\begin{equation*}
u_{\varepsilon}\rightarrow u \text{ and } v_{\varepsilon}\rightarrow v \quad \text { in } C^0_{\text {loc }}(\bar{\Omega} \times[0, \infty))\cap C^{2,1}_{\text {loc }}(\bar{\Omega} \times(0, \infty))
\end{equation*}
as $\varepsilon=\varepsilon_{j} \searrow 0,$ and that hence $(u,v)$ has the property in (\ref{Equ(1.7)}) and solves (\ref{Equ(1.2)}) classically.

\section{Global weak solutions for arbitrary $b>0$}
\hspace*{\parindent} The goal of this section is to construct the global weak solutions of (\ref{Equ(1.2)}) for any $b>0$ in the case $\gamma=2$ under certain conditions.
\subsection{Some regularity properties of approximate solutions.}
\hspace*{\parindent} In this subsection, motivated by some ideas in \cite{WAP}, we shall give a lemma which will be of great importance to establish pointwise a.e. convergence of $(u_\varepsilon)_{\varepsilon\in(0,1)}$.
\begin{lemma}\label{result4.1}
Let $\gamma=2$, and let (\ref{Equ(1.3)}), (\ref{Equ(1.12)}), (\ref{Equ(1.5)}) and (\ref{Equ(1.6)}) with some $\alpha>0$ hold. Suppose that $\beta:=\max\{1,\frac{\alpha}{2}\}$, and that $m\in \mathbb{N}$ satisfies $m>\frac{n+2}{2}$. Then one can find some $C>0$ such that
\begin{equation}\label{Equ(4.1)}
\begin{aligned}
\int_{\Omega}\left|\nabla\left(u_{\varepsilon} v_{\varepsilon}^\beta\right)\right| \leq C\left\{\int_{\Omega} \phi_{\varepsilon}\left(v_{\varepsilon}\right) \frac{\left|\nabla u_{\varepsilon}\right|^2}{u_{\varepsilon}}+\int_{\Omega} u_{\varepsilon}^2+\int_{\Omega} \frac{\left|\nabla v_{\varepsilon}\right|^2}{v_{\varepsilon}^2} +1\right\}
\end{aligned}
\end{equation}
and
\begin{equation}\label{Equ(4.2)}
\begin{aligned}
\|(u_{\varepsilon} v_{\varepsilon}^\beta)_t\|_{(W_{0}^{m,2}(\Omega))^*} \leq C\left\{\int_{\Omega} \phi_{\varepsilon}\left(v_{\varepsilon}\right) \frac{\left|\nabla u_{\varepsilon}\right|^2}{u_{\varepsilon}}+\int_{\Omega} \frac{\left|\nabla v_{\varepsilon}\right|^4}{v_{\varepsilon}^3}+\int_{\Omega} \frac{\left|\nabla v_{\varepsilon}\right|^2}{v_{\varepsilon}^2}+\int_{\Omega} |\Delta v_{\varepsilon}|^2+\int_{\Omega} u_{\varepsilon}^2+1\right\}
\end{aligned}
\end{equation}
for all $t>0$ and $\varepsilon\in(0,1)$.
\end{lemma}
\noindent{\bf{Proof.}} According to (\ref{Equ(1.5)}), (\ref{Equ(1.6)}), (\ref{Equ(1.12)}), (\ref{Equ(2.2)}) and (\ref{Equ(2.5)}), we can fix some positive constants $c_1$, $c_2$  and $c_3$ fulfilling
\begin{equation}\label{Equ(4.01)}
\begin{aligned}
c_1 v^\alpha_\varepsilon\leq\phi_\varepsilon(v_\varepsilon)\leq c_2\,\, \text{ and }\,\, |\phi_\varepsilon'(v_\varepsilon)|\leq c_{3} v^{\alpha-1}_\varepsilon \quad \text{ in } \Omega\times(0,\infty) \text{  for all  } \varepsilon\in(0,1).
\end{aligned}
\end{equation}
It follows from $\beta\geq\frac{\alpha}{2}$, (\ref{Equ(2.5)}) and (\ref{Equ(4.01)}) that
\begin{equation*}
\begin{aligned}
\int_{\Omega}\left|\nabla\left(u_{\varepsilon} v_{\varepsilon}^\beta\right)\right| \leq&\int_{\Omega} v_{\varepsilon}^\beta\left|\nabla u_{\varepsilon}\right|+\beta\int_{\Omega} u_{\varepsilon} v^{\beta-1}_{\varepsilon}\left|\nabla v_{\varepsilon}\right| \\
\leq& \frac{1}{2} \int_{\Omega} \phi_{\varepsilon}\left(v_{\varepsilon}\right) \frac{\left|\nabla u_{\varepsilon}\right|^2}{u_{\varepsilon}}+\frac{1}{2} \int_{\Omega} \frac{u_{\varepsilon} v_{\varepsilon}^{2\beta}}{\phi_{\varepsilon}\left(v_{\varepsilon}\right)}+\frac{\beta}{2} \int_{\Omega} u_{\varepsilon}^2 v_{\varepsilon}^{2\beta}+\frac{\beta}{2} \int_{\Omega} \frac{\left|\nabla v_{\varepsilon}\right|^2}{v_{\varepsilon}^2} \\
\leq& \frac{1}{2} \int_{\Omega} \phi_{\varepsilon}\left(v_{\varepsilon}\right) \frac{\left|\nabla u_{\varepsilon}\right|^2}{u_{\varepsilon}}+\frac{1}{2c_1} \int_{\Omega} u_{\varepsilon} v_{\varepsilon}^{2\beta-\alpha}+\frac{\beta}{2} \int_{\Omega} u_{\varepsilon}^2 v_{\varepsilon}^{2\beta}+\frac{\beta}{2} \int_{\Omega} \frac{\left|\nabla v_{\varepsilon}\right|^2}{v_{\varepsilon}^2} \\
\leq& \frac{1}{2} \int_{\Omega} \phi_{\varepsilon}\left(v_{\varepsilon}\right) \frac{\left|\nabla u_{\varepsilon}\right|^2}{u_{\varepsilon}}+\frac{\|v_0\|^{2\beta-\alpha}_{L^\infty(\Omega)}}{2c_1}\int_{\Omega} u_{\varepsilon}+\frac{\beta\|v_0\|^{2\beta}_{L^\infty(\Omega)}}{2} \int_{\Omega} u_{\varepsilon}^2+\frac{\beta}{2} \int_{\Omega} \frac{\left|\nabla v_{\varepsilon}\right|^2}{v_{\varepsilon}^2} \\
\end{aligned}
\end{equation*}
for all $t>0$ and $\varepsilon\in(0,1)$, this along with (\ref{Equ(2.3)}) shows (\ref{Equ(4.1)}).
Since $m>\frac{n+2}{2}$ implies the continuity of the embedding $W_{0}^{m,2}(\Omega) \hookrightarrow W^{1,\infty}(\Omega)$,
we may find $c_4>0$ with the property that $\|\varphi\|_{L^{\infty}(\Omega)} + \|\nabla \varphi\|_{L^\infty(\Omega)} \leq c_4$
for all  $\varphi \in C_0^{\infty}(\Omega)$ satisfying $\|\varphi\|_{W_0^{m,2}(\Omega)}\le1$. Fixing any such $\varphi$
and using $\beta=\max\{1,\frac{\alpha}{2}\}$, (\ref{Equ(2.5)}) and (\ref{Equ(4.01)}), we see that
\begin{equation*}
\begin{aligned}
\bigg|\int_{\Omega} \partial_t\left(u_{\varepsilon} v_{\varepsilon}^\beta\right) \varphi\bigg|=&\left|\int_{\Omega} v_{\varepsilon}^\beta\varphi u_{\varepsilon t}+\beta\int_{\Omega} u_{\varepsilon}v_{\varepsilon}^{\beta-1}\varphi v_{\varepsilon t}\right|\\
=&\bigg|-\int_{\Omega} (\phi_{\varepsilon}(v_\varepsilon)\nabla u_{\varepsilon}+u_\varepsilon\phi'_{\varepsilon}(v_\varepsilon)\nabla v_{\varepsilon})(\beta v_{\varepsilon}^{\beta-1}\varphi \nabla v_{\varepsilon}+v_{\varepsilon}^\beta\nabla\varphi)+a\int_{\Omega} u_{\varepsilon}v_{\varepsilon}^\beta\varphi\\
&-b\int_{\Omega} u^2_{\varepsilon}v_{\varepsilon}^\beta\varphi
+\beta\int_{\Omega} u_{\varepsilon} v_{\varepsilon}^{\beta-1}\varphi\Delta v_{\varepsilon}-\beta\int_{\Omega}\frac{u^2_{\varepsilon} v^\beta_{\varepsilon}}{1+\varepsilon u_{\varepsilon}} \varphi\bigg|\\
\leq& \beta c_4\int_{\Omega} \phi_{\varepsilon}\left(v_{\varepsilon}\right)v^{\beta-1}_{\varepsilon}\left|\nabla u_{\varepsilon}\right| \cdot\left|\nabla v_{\varepsilon}\right|+c_4\int_{\Omega} \phi_{\varepsilon}\left(v_{\varepsilon}\right)v_{\varepsilon}^\beta\left|\nabla u_{\varepsilon}\right|\\
&+\beta c_3c_4\int_{\Omega} u_{\varepsilon} v_{\varepsilon}^{\alpha+\beta-2}\left|\nabla v_{\varepsilon}\right|^2+c_3c_4\int_{\Omega} u_{\varepsilon} v_{\varepsilon}^{\alpha+\beta-1}\left|\nabla v_{\varepsilon}\right|\\
&+ac_4\|v_0\|^{\beta}_{L^\infty(\Omega)}\int_{\Omega} u_{\varepsilon}+bc_4\|v_0\|^{\beta}_{L^\infty(\Omega)}\int_{\Omega} u^2_{\varepsilon}\\
&+\beta c_4\int_{\Omega} u_{\varepsilon} v_{\varepsilon}^{\beta-1}\left|\Delta v_{\varepsilon}\right|+\beta c_4\int_{\Omega} u_{\varepsilon}^2 v_{\varepsilon}^\beta \\
\leq & \int_{\Omega} \phi_{\varepsilon}\left(v_{\varepsilon}\right) \frac{\left|\nabla u_{\varepsilon}\right|^2}{u_{\varepsilon}}+\frac{\beta^2c_4^2}{2}\int_{\Omega} u_{\varepsilon} \phi_{\varepsilon}\left(v_{\varepsilon}\right)v_{\varepsilon}^{2\beta-2}\left|\nabla v_{\varepsilon}\right|^2+\frac{c_4^2}{2}\int_{\Omega} u_{\varepsilon} \phi_{\varepsilon}\left(v_{\varepsilon}\right)v_{\varepsilon}^{2\beta}\\
&+\frac{1}{2}\int_{\Omega} \frac{\left|\nabla v_{\varepsilon}\right|^4}{v_{\varepsilon}^3}+\frac{\beta^2c_3^2c_4^2}{2}\int_{\Omega} u^2_{\varepsilon} v_{\varepsilon}^{2\alpha+2\beta-1}+\frac{1}{2}\int_{\Omega} \frac{\left|\nabla v_{\varepsilon}\right|^2}{v_{\varepsilon}^2}+\frac{c_3^2c_4^2}{2}\int_{\Omega} u^2_{\varepsilon} v_{\varepsilon}^{2\alpha+2\beta}\\
\end{aligned}
\end{equation*}
\begin{equation*}
\begin{aligned}
&+ac_4\|v_0\|^{\beta}_{L^\infty(\Omega)}\int_{\Omega} u_{\varepsilon}+\frac{\beta^2}{4}\int_{\Omega}\left|\Delta v_{\varepsilon}\right|^2\\
&+[(b+\beta)c_4\|v_0\|^{\beta}_{L^\infty(\Omega)}+c_4^2\|v_0\|^{2\beta-2}_{L^\infty(\Omega)}]\int_{\Omega} u^2_{\varepsilon}\\
\leq & \int_{\Omega} \phi_{\varepsilon}\left(v_{\varepsilon}\right) \frac{\left|\nabla u_{\varepsilon}\right|^2}{u_{\varepsilon}}+\int_{\Omega} \frac{\left|\nabla v_{\varepsilon}\right|^4}{v_{\varepsilon}^3}+\frac{1}{2}\int_{\Omega} \frac{\left|\nabla v_{\varepsilon}\right|^2}{v_{\varepsilon}^2}+\frac{\beta^2}{4}\int_{\Omega}\left|\Delta v_{\varepsilon}\right|^2\\
&+\bigg[\frac{\beta^4c_2^2c_4^4}{8}\|v_0\|^{4\beta-1}_{L^\infty(\Omega)}+\frac{\beta^2c_4^2 c_3^2}{2}\|v_0\|^{2\alpha+2\beta-1}_{L^\infty(\Omega)}+\frac{c_3^2c_4^2}{2}\|v_0\|^{2\alpha+2\beta}_{L^\infty(\Omega)}+c_4^2\|v_0\|^{2\beta-2}_{L^\infty(\Omega)}\\
&+(b+\beta)c_4\|v_0\|^{\beta}_{L^\infty(\Omega)}\bigg]\int_{\Omega} u^2_{\varepsilon}+\left(\frac{c_2c_4^2\|v_0\|^{2\beta}_{L^\infty(\Omega)}}{2}+ac_4\|v_0\|^{\beta}_{L^\infty(\Omega)}\right)\int_{\Omega} u_{\varepsilon}
\end{aligned}
\end{equation*}
for all $t>0$ and $\varepsilon\in(0,1)$,  which together with (\ref{Equ(2.3)}) immediately gives (\ref{Equ(4.2)}). $\hfill{} \Box$

\subsection{Weak degeneracy: $\alpha\in(0,\frac{1}{2})$. Proof of Theorem \ref{result1.2}.}
\hspace*{\parindent} When $\alpha\in(0,\frac{1}{2})$, similar to corresponding situations in (\ref{Equ(1.2)}) without the logistic source \cite{WAP}, we shall show that the uniform integrability of the corresponding first solution components will occur based on the following energy functional
$$\int_{\Omega} u_{\varepsilon} \ln u_{\varepsilon}-\delta\int_{\Omega} u_{\varepsilon} \phi_{\varepsilon}\left(v_{\varepsilon}\right),$$
where $\delta$ is some positive constant and will be determined later. The following inequality of great importance has been recorded in \cite[Lemma 3.3]{WAP}.

\begin{lemma}\label{result4.2}
Let (\ref{Equ(1.12)}) and (\ref{Equ(1.13)}) be satisfied with some $\alpha \in(0,1)$ and $s_0>0$. Then there exist $\delta>0$ and $\varepsilon_{0} \in(0,1)$ with the property that for any $s_{\star}>0$, we can find $C\left(s_{\star}\right)>0$ fulfilling
$$
\frac{\left\{\delta\phi_{\varepsilon}(s) \phi_{\varepsilon}^{\prime}(s)+(\delta-1) \phi_{\varepsilon}^{\prime}(s)\right\}^2}{2 \phi_{\varepsilon}(s)}+\delta \phi_{\varepsilon}^{\prime 2}(s)+\delta \phi_{\varepsilon}^{\prime \prime}(s) \leq \frac{C\left(s_{\star}\right)}{s} \quad \text { for all } s \in\left[0, s_{\star}\right] \text { and } \varepsilon \in\left(0, \varepsilon_{0}\right).
$$
\end{lemma}

\begin{lemma}\label{result4.3}
Let $\gamma=2$, and assume that (\ref{Equ(1.3)}) and (\ref{Equ(1.12)}) hold. Then we have
\begin{equation*}
\begin{aligned}
\frac{d}{d t} \int_{\Omega} u_{\varepsilon} \phi_{\varepsilon}\left(v_{\varepsilon}\right)=&-\int_{\Omega}\left\{\phi_{\varepsilon}\left(v_{\varepsilon}\right) \phi_{\varepsilon}^{\prime}\left(v_{\varepsilon}\right)+\phi_{\varepsilon}^{\prime}\left(v_{\varepsilon}\right)\right\} \nabla u_{\varepsilon} \cdot \nabla v_{\varepsilon}-\int_{\Omega} u_{\varepsilon}\left\{\phi_{\varepsilon}^{\prime 2}\left(v_{\varepsilon}\right)+\phi_{\varepsilon}^{\prime \prime}\left(v_{\varepsilon}\right)\right\}\left|\nabla v_{\varepsilon}\right|^2 \\
&-\int_{\Omega} \frac{u_{\varepsilon}^2 v_{\varepsilon}}{1+\varepsilon u_{\varepsilon}} \phi_{\varepsilon}^{\prime}\left(v_{\varepsilon}\right)+a\int_{\Omega} u_{\varepsilon} \phi_{\varepsilon}\left(v_{\varepsilon}\right)-b\int_{\Omega} u^2_{\varepsilon} \phi_{\varepsilon}\left(v_{\varepsilon}\right)
\end{aligned}
\end{equation*}
for all $t>0$ and $\varepsilon \in(0,1)$.
\end{lemma}
\noindent{\bf{Proof.}} We take a similar procedure as the proof of \cite[Lemma 3.2]{WAP} to achieve this assertion. $\hfill{} \Box$

\begin{lemma}\label{result4.4}
Under the assumptions of Theorem \ref{result1.2}, there exist $\varepsilon_{0} \in(0,1)$ and $C>0$ satisfying
\begin{equation}\label{Equ(4.3)}
\int_\Omega u_\varepsilon(\cdot,t)\ln u_\varepsilon(\cdot,t)\leq C \quad \text{for all } t>0 \text{ and }\varepsilon\in(0,\varepsilon_{0}),
\end{equation}
and
\begin{equation}\label{Equ(4.4)}
\int_{t}^{t+1}\int_\Omega \phi_\varepsilon(v_\varepsilon)\frac{|\nabla u_\varepsilon|^2}{u_\varepsilon}+\int_{t}^{t+1}\int_\Omega u^2_\varepsilon\ln u_\varepsilon\leq C \quad \text{for all } t>0 \text{ and }\varepsilon\in(0,\varepsilon_{0}).
\end{equation}
\end{lemma}
\noindent{\bf{Proof.}} Suppose that $\delta>0$ and $\varepsilon_{0}\in(0,1)$ are taken from Lemma \ref{result4.2}. Lemma \ref{result4.2} can be applied to $s_{\star}:=\|v_0\|_{L^\infty(\Omega)}$ to show that there exists $c_1>0$ fulfilling
\begin{equation}\label{Equ(4.5)}
\frac{\left\{\delta\phi_{\varepsilon}(v_\varepsilon) \phi_{\varepsilon}^{\prime}(v_\varepsilon)+(\delta-1) \phi_{\varepsilon}^{\prime}(v_\varepsilon)\right\}^2}{2 \phi_{\varepsilon}(v_\varepsilon)}+\delta \phi_{\varepsilon}^{\prime 2}(v_\varepsilon)+\delta \phi_{\varepsilon}^{\prime \prime}(v_\varepsilon) \leq \frac{c_1}{v_\varepsilon}\,\,\text { in } \Omega\times(0,\infty) \text { for all } \varepsilon \in\left(0, \varepsilon_{0}\right).
\end{equation}
Note that $(a+1)\int_{\Omega} s \ln s-\frac{b}{2}\int_{\Omega} s^2 \ln s+a\int_{\Omega} s-b\int_{\Omega} s^2$ for all $s>0$ is bounded by some constant $c_2>0$. Testing the first equation of (\ref{Equ(2.1)}) by $1+\ln u_\varepsilon$ and then integrating by parts, we show that
\begin{equation}\label{Equ(4.6)}
\begin{aligned}
\frac{d}{d t}&\int_{\Omega} u_{\varepsilon} \ln u_{\varepsilon}+\int_{\Omega} u_{\varepsilon} \ln u_{\varepsilon}+\int_{\Omega} \phi_\varepsilon\left(v_{\varepsilon}\right) \frac{\left|\nabla u_{\varepsilon}\right|^{2}}{u_{\varepsilon}}
+\frac{b}{2}\int_{\Omega} u^2_{\varepsilon} \ln u_{\varepsilon}\\
=&-\int_{\Omega} \phi_\varepsilon^{\prime}\left(v_{\varepsilon}\right) \nabla u_{\varepsilon} \cdot \nabla v_{\varepsilon}+(a+1)\int_{\Omega} u_{\varepsilon} \ln u_{\varepsilon}-\frac{b}{2}\int_{\Omega} u^2_{\varepsilon} \ln u_{\varepsilon}+a\int_{\Omega} u_{\varepsilon}-b\int_{\Omega} u^2_{\varepsilon}\\
\leq&-\int_{\Omega} \phi_\varepsilon^{\prime}\left(v_{\varepsilon}\right) \nabla u_{\varepsilon} \cdot \nabla v_{\varepsilon}+c_{2}
\end{aligned}
\end{equation}
for all $ t>0 $ and $\varepsilon \in(0,1)$. By (\ref{Equ(1.6)}), (\ref{Equ(1.12)}), (\ref{Equ(2.2)}) and (\ref{Equ(2.5)}), one can find some $c_3>0$ and $c_4>0$  such that
$\phi_\varepsilon(v_\varepsilon)\leq c_3$ and $ |\phi_\varepsilon'(v_\varepsilon)|\leq c_{4} v^{\alpha-1}_\varepsilon$  in $\Omega\times(0,\infty)$  for all $\varepsilon\in(0,1)$. Then we utilize Lemma \ref{result4.3}, (\ref{Equ(4.5)}), (\ref{Equ(4.6)}), (\ref{Equ(2.5)}) and Young's inequality to obtain
\begin{equation*}
\begin{aligned}
\frac{d}{d t}&\left\{\int_{\Omega} u_{\varepsilon} \ln u_{\varepsilon}-\delta\int_{\Omega} u_{\varepsilon} \phi_{\varepsilon}\left(v_{\varepsilon}\right)\right\}+\int_{\Omega} u_{\varepsilon} \ln u_{\varepsilon}-\delta\int_{\Omega} u_{\varepsilon} \phi_{\varepsilon}\left(v_{\varepsilon}\right)\\
&+\int_{\Omega} \phi_{\varepsilon}\left(v_{\varepsilon}\right) \frac{\left|\nabla u_{\varepsilon}\right|^2}{u_{\varepsilon}}+\frac{b}{2}\int_{\Omega} u^2_{\varepsilon} \ln u_{\varepsilon}\\
\leq&\int_{\Omega}\left\{\delta\phi_{\varepsilon}\left(v_{\varepsilon}\right) \phi_{\varepsilon}^{\prime}\left(v_{\varepsilon}\right)+(\delta-1)\phi_{\varepsilon}^{\prime}\left(v_{\varepsilon}\right)\right\} \nabla u_{\varepsilon} \cdot \nabla v_{\varepsilon}+\delta\int_{\Omega} u_{\varepsilon}\left\{\phi_{\varepsilon}^{\prime 2}\left(v_{\varepsilon}\right)+\phi_{\varepsilon}^{\prime \prime}\left(v_{\varepsilon}\right)\right\}\left|\nabla v_{\varepsilon}\right|^2 \\
&+\delta\int_{\Omega} \frac{u_{\varepsilon}^2 v_{\varepsilon}}{1+\varepsilon u_{\varepsilon}} \phi_{\varepsilon}^{\prime}\left(v_{\varepsilon}\right)-\delta(a+1)\int_{\Omega} u_{\varepsilon} \phi_{\varepsilon}\left(v_{\varepsilon}\right)+b\delta\int_{\Omega} u^2_{\varepsilon} \phi_{\varepsilon}\left(v_{\varepsilon}\right)+c_2\\
\leq&\frac{1}{2}\int_{\Omega} \phi_{\varepsilon}\left(v_{\varepsilon}\right) \frac{\left|\nabla u_{\varepsilon}\right|^2}{u_{\varepsilon}}+\int_{\Omega}\left\{\frac{(\delta\phi_{\varepsilon}(v_{\varepsilon}) \phi_{\varepsilon}^{\prime}(v_{\varepsilon})+(\delta-1)\phi_{\varepsilon}^{\prime}(v_{\varepsilon}) )^2}{2\phi_{\varepsilon}(v_{\varepsilon})}+\delta\phi_{\varepsilon}^{\prime 2}(v_{\varepsilon})+\delta\phi_{\varepsilon}^{\prime \prime}(v_{\varepsilon})\right\}u_{\varepsilon}|\nabla v_{\varepsilon}|^2\\
&+\delta\int_{\Omega} \frac{u_{\varepsilon}^2 v_{\varepsilon}}{1+\varepsilon u_{\varepsilon}} \phi_{\varepsilon}^{\prime}\left(v_{\varepsilon}\right)-\delta(a+1)\int_{\Omega} u_{\varepsilon} \phi_{\varepsilon}\left(v_{\varepsilon}\right)+b\delta\int_{\Omega} u^2_{\varepsilon} \phi_{\varepsilon}\left(v_{\varepsilon}\right)+c_2\\
\end{aligned}
\end{equation*}
\begin{equation}\label{Equ(4.7)}
\begin{aligned}
\leq&\frac{1}{2}\int_{\Omega} \phi_{\varepsilon}\left(v_{\varepsilon}\right) \frac{\left|\nabla u_{\varepsilon}\right|^2}{u_{\varepsilon}}+c_1\int_{\Omega}\frac{u_{\varepsilon}}{v_{\varepsilon}}|\nabla v_{\varepsilon}|^2+\delta c_4\|v_0\|^{\alpha}_{L^\infty(\Omega)}\int_{\Omega} u_{\varepsilon}^2+b\delta c_3\int_{\Omega} u^2_{\varepsilon}+c_2\\
\leq&\frac{1}{2}\int_{\Omega} \phi_{\varepsilon}\left(v_{\varepsilon}\right) \frac{\left|\nabla u_{\varepsilon}\right|^2}{u_{\varepsilon}}+\frac{c^2_1}{4}\int_{\Omega}\frac{|\nabla v_{\varepsilon}|^4}{v^3_{\varepsilon}}+\left\{\|v_0\|_{L^\infty(\Omega)}+\delta c_4\|v_0\|^{\alpha}_{L^\infty(\Omega)}+b\delta c_3\right\}\int_{\Omega} u_{\varepsilon}^2+c_2
\end{aligned}
\end{equation}
for all $ t>0 $ and $\varepsilon \in(0,\varepsilon_{0})$. Hence, using (\ref{Equ(4.7)}), (\ref{Equ(1.12)}), (\ref{Equ(2.3)}), (\ref{Equ(2.4)}), (\ref{Equ(2.5)}), (\ref{Equ(2.15)}) and \cite[Lemma 3.4]{SSW} we obtain (\ref{Equ(4.3)}). Integrating (\ref{Equ(4.7)}), we infer (\ref{Equ(4.4)}) from  (\ref{Equ(1.12)}), (\ref{Equ(2.3)}) and (\ref{Equ(2.5)}).
Hence, we complete this proof. $\hfill{} \Box$

\begin{corollary}\label{result4.01}
Let the assumptions of Theorem \ref{result1.2} hold. Suppose that $m\in \mathbb{N}$ satisfies $m>\frac{n+2}{2}$. Then we can find $\varepsilon_{0}\in(0,1)$ such that for each $T>0$ there exists $C(T)>0$ fulfilling
\begin{equation*}
\int_{0}^T\int_{\Omega}\left|\nabla\left(u_{\varepsilon} v_{\varepsilon}\right)\right| \leq C(T)\quad \text{for all } t\in(0,T) \text{ and }\varepsilon\in(0,\varepsilon_{0})
\end{equation*}
and
\begin{equation*}
\int_{0}^T\|(u_{\varepsilon} v_{\varepsilon})_t\|_{(W_{0}^{m,2}(\Omega))^*} \leq C(T)\quad \text{for all } t\in(0,T) \text{ and }\varepsilon\in(0,\varepsilon_{0}).
\end{equation*}
\end{corollary}
\noindent{\bf{Proof.}} Both assertions are a consequence of Lemmata \ref{result4.1}, \ref{result2.6}, (\ref{Equ(2.4)}), (\ref{Equ(2.5)}) and (\ref{Equ(4.4)}).  $\hfill{} \Box$

Based on (\ref{Equ(2.4)}) and Lemmata \ref{result2.6} and \ref{result4.4}, we can make sure that $\nabla(u_\varepsilon\phi_\varepsilon(v_\varepsilon))$ has the following integrability feature in reflexive Lebesgue spaces.
\begin{lemma}\label{result4.5}
Under the assumptions of Theorem \ref{result1.2}, given any $T>0$, there exist $\varepsilon_{0} \in(0,1)$ and $C(T)>0$ such that
\begin{equation*}
\int_{0}^{T}\int_\Omega|\nabla(u_\varepsilon\phi_\varepsilon(v_\varepsilon))|^{p(\alpha)}\leq C(T) \quad \text{for all } t\in(0,T) \text{ and } \varepsilon \in(0,\varepsilon_{0}),
\end{equation*}
where
\begin{equation*}
p(\alpha):=
\begin{cases}\frac{1}{1-\alpha} & \text { if } \alpha \in\left(0, \frac{1}{4}\right),\\
\frac{4}{3} & \text { if } \alpha \in\left[\frac{1}{4}, \frac{1}{2}\right).
\end{cases}
\end{equation*}
\end{lemma}
\noindent{\bf{Proof.}} Let $p=p(\alpha)$. It follows from (\ref{Equ(1.6)}), (\ref{Equ(2.2)}) and (\ref{Equ(2.5)}) that there is $c_1>0$ satisfying
$|\phi_\varepsilon'(v_\varepsilon)|\leq c_1 v^{\alpha-1}_\varepsilon$ in $\Omega\times(0,\infty)$  for all $\varepsilon\in(0,1)$. Straightforward calculations show that
\begin{equation}\label{Equ(4.8)}
\begin{aligned}
\int_\Omega|\nabla(u_\varepsilon\phi_\varepsilon(v_\varepsilon))|^p
\leq& 2^{p-1}\int_{\Omega}\phi^p_\varepsilon(v_\varepsilon)|\nabla u_{\varepsilon}|^p+2^{p-1}c_1^{p}\int_{\Omega} u^p_{\varepsilon}v^{p(\alpha-1)}_\varepsilon |\nabla v_{\varepsilon}|^p\\
=:& 2^{p-1}I_\varepsilon(t)+2^{p-1}c_1^{p}J_\varepsilon(t)
\end{aligned}
\end{equation}
for all $t>0$ and $\varepsilon \in(0,1)$. Since $p\leq\frac{4}{3}$, we can readily obtain $\frac{p}{2-p}\leq2$. Hence, in light of Young's inequality, we infer that
\begin{equation}\label{Equ(4.9)}
\begin{aligned}
I_\varepsilon(t)=&\int_{\Omega} \left(\phi_{\varepsilon}\left(v_{\varepsilon}\right) \frac{\left|\nabla u_{\varepsilon}\right|^2}{u_{\varepsilon}}\right)^\frac{p}{2} (u_{\varepsilon}\phi_\varepsilon(v_\varepsilon))^\frac{p}{2}\\
\leq&\int_{\Omega} \phi_{\varepsilon}\left(v_{\varepsilon}\right) \frac{\left|\nabla u_{\varepsilon}\right|^2}{u_{\varepsilon}}+\int_{\Omega} (u_{\varepsilon}\phi_\varepsilon(v_\varepsilon))^\frac{p}{2-p}\\
\leq&\int_{\Omega} \phi_{\varepsilon}\left(v_{\varepsilon}\right) \frac{\left|\nabla u_{\varepsilon}\right|^2}{u_{\varepsilon}}+c_2^2\int_{\Omega} u_{\varepsilon}^2
+|\Omega|^\frac{4-3p}{4-2p}
\end{aligned}
\end{equation}
for all $t>0$ and $\varepsilon \in(0,1)$ with $c_2:=\sup_{\varepsilon \in(0,1)}\|\phi_\varepsilon(v_\varepsilon)\|_{L^\infty(\Omega\times(0,\infty))}$ being finite due to (\ref{Equ(1.12)}), (\ref{Equ(2.2)}) and (\ref{Equ(2.5)}). For $J_\varepsilon(t)$, we have two cases to proceed as below: In the case $\alpha\in[\frac{1}{4},\frac{1}{2})$, we employ (\ref{Equ(2.5)}) and
Young's inequality to show that
\begin{equation}\label{Equ(4.10)}
\begin{aligned}
J_\varepsilon(t)=\int_{\Omega} \left(\frac{\left|\nabla v_{\varepsilon}\right|^4}{v^3_{\varepsilon}}\right)^\frac{1}{3} u^\frac{4}{3}_{\varepsilon}v_\varepsilon^\frac{4\alpha-1}{3}
\leq&\int_{\Omega} \frac{\left|\nabla v_{\varepsilon}\right|^4}{v^3_{\varepsilon}}+\int_{\Omega}u^2_{\varepsilon}v_\varepsilon^\frac{4\alpha-1}{2}\\
\leq&\int_{\Omega} \frac{\left|\nabla v_{\varepsilon}\right|^4}{v^3_{\varepsilon}}+\|v_0\|^\frac{4\alpha-1}{2}_{L^\infty(\Omega)}\int_{\Omega}u^2_{\varepsilon}
\end{aligned}
\end{equation}
for all $t>0$ and $\varepsilon \in(0,1)$, whereas in the case $\alpha\in(0,\frac{1}{4})$, an application of $p=\frac{1}{1-\alpha}$ and Young's inequality entails that
\begin{equation}\label{Equ(4.11)}
\begin{aligned}
J_\varepsilon(t)=&\int_{\Omega} u^\frac{1}{1-\alpha}_{\varepsilon}v^{-1}_\varepsilon |\nabla v_{\varepsilon}|^\frac{1}{1-\alpha}\\
\leq&\int_{\Omega}u^2_{\varepsilon}+\int_{\Omega}v_\varepsilon^\frac{2(\alpha-1)}{1-2\alpha}|\nabla v_{\varepsilon}|^\frac{2}{1-2\alpha}\\
=&\int_{\Omega}u^2_{\varepsilon}+\int_{\Omega}\left(\frac{|\nabla v_{\varepsilon}|^2}{v_\varepsilon^2}\right)^\frac{1-4\alpha}{1-2\alpha}\left(\frac{|\nabla v_{\varepsilon}|^4}{v_\varepsilon^3}\right)^\frac{2\alpha}{1-2\alpha}\\
\leq&\int_{\Omega}u^2_{\varepsilon}+\int_{\Omega}\frac{|\nabla v_{\varepsilon}|^2}{v_\varepsilon^2}+\int_{\Omega}\frac{|\nabla v_{\varepsilon}|^4}{v_\varepsilon^3}\quad \text{for all } t>0 \text{ and } \varepsilon \in(0,1).
\end{aligned}
\end{equation}
Combining (\ref{Equ(4.8)})-(\ref{Equ(4.11)}), (\ref{Equ(2.4)}), (\ref{Equ(2.15)}), (\ref{Equ(2.17)}) and (\ref{Equ(4.4)}), we can achieve the intended conclusion.  $\hfill{} \Box$

With the help of a priori estimates, we can extract appropriately subsequences from $(u_{\varepsilon}, v_{\varepsilon})$ such that it is convergent.
\begin{lemma}\label{result4.6}
Let the assumptions of Theorem \ref{result1.2} hold. Then there exist a sequence $\left(\varepsilon_{j}\right)_{j \in \mathbb{N}} \subset(0,1)$ and functions $u$ and $v$ which fulfill (\ref{Equ(1.14)}) with $u \geq 0$ and $v>0$ a.e. in $\Omega \times(0, \infty)$, such that $\varepsilon_{j} \searrow 0$ as $j \rightarrow \infty$ and that
\begin{align}
&u_{\varepsilon}\rightarrow u \quad \text { in } L^2_{\text {loc }}(\bar{\Omega} \times[0, \infty)) \text { and } \text { a.e. in } \Omega\times (0, \infty),\label{Equ(4.13)}\\
&v_{\varepsilon} \rightarrow v \quad \text { in } L^2_{\text {loc }}([0, \infty), W^{1,2}(\Omega))\cap C^0_{\text {loc }}([0, \infty), L^p(\Omega))\nonumber\\
&\quad\quad\quad\quad\text { for all } p\in\left[1,\frac{2n}{(n-2)_+}\right)\text { and } \text { a.e. in } \Omega\times (0, \infty),\label{Equ(4.14)}\\
&\nabla v_{\varepsilon} \rightharpoonup \nabla v \quad \text { in } L^4_{\text {loc }}(\bar{\Omega} \times[0, \infty)) ,\label{Equ(4.15)}\\
&v_{\varepsilon} \stackrel{\star}{\rightharpoonup} v \quad \text { in } L^\infty(\Omega\times(0, \infty)), \label{Equ(4.16)}\\
&\frac{u_\varepsilon v_\varepsilon}{1+\varepsilon u_\varepsilon} \rightarrow uv \text { in } L_{\text {loc }}^1(\bar{\Omega} \times[0, \infty)),\label{Equ(4.17)}\\
&\nabla(\phi(v_{\varepsilon})u_{\varepsilon}) \rightharpoonup \nabla(\phi(v)u) \quad \text { in } L^{p(\alpha)}_{\text {loc }}(\bar{\Omega} \times[0, \infty))\label{Equ(4.18)}
\end{align}
as $\varepsilon=\varepsilon_{j} \searrow 0,$ where $p(\alpha)$ is taken from Lemma \ref{result4.5}, and that $(u,v)$ forms a global weak solution of (\ref{Equ(1.8)})-(\ref{Equ(1.11)}) in the sense of Definition \ref{result0.1}.
\end{lemma}
\noindent{\bf{Proof.}} Let $ T>0$. Using Corollary \ref{result4.01} along with (\ref{Equ(2.3)}) and (\ref{Equ(2.5)}), and letting $m \in \mathbb{N}$ such that $m>\frac{n+2}{2}$, we infer the existence of $\varepsilon_{0}\in (0,1)$ fulfilling
$$(u_{\varepsilon}v_{\varepsilon})_{\varepsilon \in(0,\varepsilon_{0})} \text{ is bounded in } L_{loc}^1\left([0, \infty) ; W^{1, 1}(\Omega)\right)$$
and $$((u_{\varepsilon}v_{\varepsilon})_t)_{\varepsilon \in(0,\varepsilon_{0})} \text{ is bounded in } L^1_{loc}\left([0, \infty) ;\left(W_{0}^{m,2}(\Omega)\right)^{\star}\right).$$
An application of an Aubin-Lions type lemma (\cite[Corollary 4]{SC}) implies that
\begin{equation*}
(u_{\varepsilon}v_{\varepsilon})_{\varepsilon \in(0,\varepsilon_{0})} \text{ is relatively compact with respect to the strong topology in } L_{loc}^1\left(\overline{\Omega}\times[0, \infty)\right),
\end{equation*}
which immediately gives the existence of a nonnegative function $z\in L_{loc}^1\left(\overline{\Omega}\times[0, \infty)\right)$ and $\left(\varepsilon_{j}\right)_{j \in \mathbb{N}} \subset(0,\varepsilon_{0})$  such that $\varepsilon_{j} \searrow 0$ as $j\rightarrow\infty$, and that
\begin{equation}\label{Equ(4.19)}
u_{\varepsilon}v_{\varepsilon}\rightarrow z \quad a.e. \text{ in } \Omega\times(0, \infty) \quad\text{ as } \varepsilon=\varepsilon_{j} \searrow0.
\end{equation}
Lemma \ref{result4.4} shows that
\begin{equation}\label{Equ(4.20)}
\begin{aligned}
(u_{\varepsilon})_{\varepsilon \in(0,\varepsilon_{0})} \text { and } (u^2_{\varepsilon})_{\varepsilon \in(0,\varepsilon_{0})}\text{ are uniformly integrable over } \Omega\times(0,T).
\end{aligned}
\end{equation}
In view of the Dunford-Pettis theorem, $(u_{\varepsilon})_{\varepsilon \in(0,\varepsilon_{0})}$ and $(u^2_{\varepsilon})_{\varepsilon \in(0,\varepsilon_{0})}$ are relatively compact with respect to the weak topology in $L^1(\Omega \times(0, T))$.

We utilize (\ref{Equ(2.5)}), (\ref{Equ(2.12)})-(\ref{Equ(2.14)}) to infer that $(v_{\varepsilon})_{\varepsilon \in(0,1)}$ is bounded in $L^\infty\left((0, \infty) ; W^{1,2}(\Omega)\right)\cap L_{loc}^{2}\left([0, \infty) ; W^{2,2}(\Omega)\right)$ and that $(v_{\varepsilon t})_{\varepsilon \in(0,1)}$ is bounded in $L_{loc}^{2}([0, \infty) ;L^2(\Omega))$. Therefore, an Aubin-Lions lemma guarantees (\ref{Equ(4.14)})-(\ref{Equ(4.16)}) with a nonnegative function $v\in L^\infty(\Omega\times(0, \infty))\cap L^\infty\left((0, \infty) ; W^{1,2}(\Omega)\right)\cap L_{loc}^{2}\left([0, \infty) ; W^{2,2}(\Omega)\right)$. Hence, (\ref{Equ(2.16)}) in conjunction with (\ref{Equ(4.14)}) and Fatou's lemma shows that $v>0$ a.e. in $\Omega\times(0,\infty)$. Letting $u:=\frac{z}{v}$ a.e. in $\Omega\times(0,\infty)$, and using (\ref{Equ(4.19)}) and (\ref{Equ(4.14)}), we infer that $u_\varepsilon\rightarrow u$ a.e. in $\Omega\times(0,\infty)$ as $\varepsilon_{j} \searrow 0$, which together with (\ref{Equ(2.3)}) and (\ref{Equ(2.4)}) shows the nonnegative function $u \in L^\infty((0,\infty); L^1(\Omega))\cap L_{loc}^2([0,\infty); L^2(\Omega))$. Therefore, (\ref{Equ(4.20)}) entails that we can extract subsequence such that $u^2_\varepsilon\rightharpoonup u^2$ in $ L^1\left(\Omega\times(0, T)\right)$, which immediately shows $u^2_\varepsilon\rightarrow u^2$ in $ L^1\left(\Omega\times(0, T)\right)$. Moreover, along with (\ref{Equ(2.4)}) and a further subsequence, we have $u_\varepsilon\rightharpoonup u$ in $L_{loc}^{2}([0, \infty); L^2(\Omega))$ and then (\ref{Equ(4.13)}).

We note that (\ref{Equ(2.5)}) and (\ref{Equ(4.20)}) entail that $(\frac{u_\varepsilon v_\varepsilon}{1+\varepsilon u_\varepsilon})_{\varepsilon \in(0,\varepsilon_{0})}$  is uniformly integrable over $\Omega\times(0,T)$. Together with (\ref{Equ(4.13)}), (\ref{Equ(4.14)}) and again the Vitali convergence theorem we obtain (\ref{Equ(4.17)}). Furthermore, (\ref{Equ(1.12)}), (\ref{Equ(2.5)}) and (\ref{Equ(4.20)}) show that $(u_\varepsilon\phi(v_\varepsilon))_{\varepsilon \in(0,\varepsilon_{0})}$  is uniformly integrable over $\Omega\times(0,T)$, by (\ref{Equ(4.13)}) and (\ref{Equ(4.14)}) we see that $u_\varepsilon\phi(v_\varepsilon) \rightarrow u\phi(v)$ in $L_{loc}^1(\bar{\Omega} \times[0, \infty))$, which along with Lemma \ref{result4.5} immediately implies (\ref{Equ(4.18)}). Because these convergence properties are sufficient to pass to the limit in each integral making up a weak formulation of (\ref{Equ(2.1)}), so $(u,v)$ is a global weak solution of (\ref{Equ(1.2)}). $\hfill{} \Box$

We are now in the position to prove Theorem \ref{result1.2}.\\
{\bf Proof of Theorem \ref{result1.2}.} Theorem \ref{result1.2} is a consequence of Lemma \ref{result4.6}. $\hfill{} \Box$

\subsection{Strong degeneracy: $\alpha\geq\frac{1}{2}$. Proof of Theorem \ref{result1.3}.}

\hspace*{\parindent}In this subsection, we will construct global weak solutions in the case $\alpha\geq\frac{1}{2}$. To this end, we first derive an entropy-like inequality to yield the following fundamental estimates.

\begin{lemma}\label{result4.7}
Let $\gamma=2$ and $\alpha\geq\frac{1}{2}$, and suppose that (\ref{Equ(1.3)}), (\ref{Equ(1.5)}), (\ref{Equ(1.6)}) and (\ref{Equ(1.12)}) hold. Then there exists $C>0$ such that
\begin{equation}\label{Equ(4.21)}
\int_\Omega u_\varepsilon(\cdot,t)\ln u_\varepsilon(\cdot,t)\leq C \quad \text{for all } t>0 \text{ and }\varepsilon\in(0,1)
\end{equation}
and
\begin{equation}\label{Equ(4.22)}
\int_{t}^{t+1}\int_\Omega \phi_\varepsilon(v_\varepsilon)\frac{|\nabla u_\varepsilon|^2}{u_\varepsilon}+\int_{t}^{t+1}\int_\Omega u^2_\varepsilon\ln u_\varepsilon\leq C \quad \text{for all } t>0 \text{ and }\varepsilon\in(0,1).
\end{equation}
\end{lemma}
\noindent{\bf{Proof.}} For all $s>0$, it can readily be verified that $(a+1)\int_{\Omega} s \ln s-\frac{b}{2}\int_{\Omega} s^\gamma \ln s+a\int_{\Omega} s-b\int_{\Omega} s^2\leq c_1$ with some $c_1>0$.  Once more relying on (\ref{Equ(1.5)}), (\ref{Equ(1.6)}),  (\ref{Equ(2.2)}) and (\ref{Equ(2.5)}), there exist $c_2>0$  and $c_3>0$  such that
$c_2 v^\alpha_\varepsilon\leq\phi_\varepsilon(v_\varepsilon)$ and $ |\phi_\varepsilon'(v_\varepsilon)|\leq c_{3} v^{\alpha-1}_\varepsilon$ in $ \Omega\times(0,\infty)$   for all  $\varepsilon\in(0,1).$ Then by a direct calculation, we see that
\begin{equation}\label{Equ(4.23)}
\begin{aligned}
\frac{d}{d t}&\int_{\Omega} u_{\varepsilon} \ln u_{\varepsilon}+\int_{\Omega} u_{\varepsilon} \ln u_{\varepsilon}+\int_{\Omega} \phi_\varepsilon\left(v_{\varepsilon}\right) \frac{\left|\nabla u_{\varepsilon}\right|^{2}}{u_{\varepsilon}}
+\frac{b}{2}\int_{\Omega} u^2_{\varepsilon} \ln u_{\varepsilon}\\
=&-\int_{\Omega} \phi_\varepsilon^{\prime}\left(v_{\varepsilon}\right) \nabla u_{\varepsilon} \cdot \nabla v_{\varepsilon}+(a+1)\int_{\Omega} u_{\varepsilon} \ln u_{\varepsilon}-\frac{b}{2}\int_{\Omega} u^2_{\varepsilon} \ln u_{\varepsilon}+a\int_{\Omega} u_{\varepsilon}-b\int_{\Omega} u^2_{\varepsilon}\\
\leq&\frac{1}{2}\int_{\Omega} \phi_\varepsilon\left(v_{\varepsilon}\right) \frac{\left|\nabla u_{\varepsilon}\right|^{2}}{u_{\varepsilon}}+\frac{1}{2}\int_{\Omega} \frac{|\phi'_\varepsilon\left(v_{\varepsilon}\right)|^2}{\phi_\varepsilon\left(v_{\varepsilon}\right)}u_{\varepsilon}|\nabla v_{\varepsilon}|^{2}+c_{1}\\
\leq&\frac{1}{2}\int_{\Omega} \phi_\varepsilon\left(v_{\varepsilon}\right) \frac{\left|\nabla u_{\varepsilon}\right|^{2}}{u_{\varepsilon}}+\frac{c_3^2}{2c_2}\int_{\Omega}u_{\varepsilon}v^{\alpha-2}_{\varepsilon}|\nabla v_{\varepsilon}|^{2}+c_{1}\\
\leq&\frac{1}{2}\int_{\Omega} \phi_\varepsilon\left(v_{\varepsilon}\right) \frac{\left|\nabla u_{\varepsilon}\right|^{2}}{u_{\varepsilon}}+\frac{c_3^2}{4c_2}\int_{\Omega}u^2_{\varepsilon}+\frac{c_3^2}{4c_2}\int_{\Omega}|\nabla v_{\varepsilon}|^{4}v^{2\alpha-4}_{\varepsilon}+c_{1}\\
\leq&\frac{1}{2}\int_{\Omega} \phi_\varepsilon\left(v_{\varepsilon}\right) \frac{\left|\nabla u_{\varepsilon}\right|^{2}}{u_{\varepsilon}}+\frac{c_3^2}{4c_2}\int_{\Omega}u^2_{\varepsilon}+\frac{c_3^2\|v_0\|^{2\alpha-1}_{L^\infty(\Omega)}}{4c_2}\int_{\Omega}\frac{|\nabla v_{\varepsilon}|^{4}}{v^3_{\varepsilon}}+c_{1}
\end{aligned}
\end{equation}
for all $ t>0 $ and $\varepsilon \in(0,1)$, because the case $\alpha\geq\frac{1}{2}$. Using (\ref{Equ(2.4)}) and (\ref{Equ(2.15)}) we derive (\ref{Equ(4.21)}), then integrating (\ref{Equ(4.23)}) we have (\ref{Equ(4.22)}). Hence, we complete this proof. $\hfill{} \Box$

Now the following estimates can be readily proved.
\begin{corollary}\label{result4.02}
Let $\gamma=2$ and $\alpha\geq\frac{1}{2}$, and let (\ref{Equ(1.3)}), (\ref{Equ(1.5)}), (\ref{Equ(1.6)}) and (\ref{Equ(1.12)}) hold. Suppose that $\beta:=\max\{1,\frac{\alpha}{2}\}$, and that $m\in \mathbb{N}$ fulfills $m>\frac{n+2}{2}$. Then for each $T>0$ there exists $C(T)>0$ satisfying
\begin{equation*}
\int_{0}^T\int_{\Omega}\left|\nabla\left(u_{\varepsilon} v^\beta_{\varepsilon}\right)\right| \leq C(T)\quad \text{for all } t\in(0,T) \text{ and }\varepsilon\in(0,1)
\end{equation*}
and
\begin{equation*}
\int_{0}^T\|(u_{\varepsilon} v^\beta_{\varepsilon})_t\|_{(W_{0}^{m,2}(\Omega))^*} \leq C(T)\quad \text{for all } t\in(0,T) \text{ and }\varepsilon\in(0,1).
\end{equation*}
\end{corollary}
\noindent{\bf{Proof.}} These directly result from (\ref{Equ(2.4)}), (\ref{Equ(2.5)}), (\ref{Equ(4.22)}), Lemmata \ref{result4.1} and \ref{result2.6}.  $\hfill{} \Box$

A combination of (\ref{Equ(2.4)}), (\ref{Equ(2.15)}) and (\ref{Equ(4.22)}) asserts a regularity feature of $\nabla(u_\varepsilon\phi_\varepsilon(v_\varepsilon))$.
\begin{lemma}\label{result4.8}
Let $\gamma=2$, and suppose that (\ref{Equ(1.3)}), (\ref{Equ(1.5)}), (\ref{Equ(1.6)}) and (\ref{Equ(1.12)}) with $\alpha\geq\frac{1}{2}$ hold. Then we can find some $C>0$ such that
\begin{equation*}
\int_{t}^{t+1}\int_\Omega|\nabla(u_\varepsilon\phi_\varepsilon(v_\varepsilon))|^\frac{4}{3}\leq C \quad \text{for all } t>0 \text{ and } \varepsilon \in(0,1).
\end{equation*}
\end{lemma}
\noindent{\bf{Proof.}} According to (\ref{Equ(1.6)}), (\ref{Equ(1.12)}), (\ref{Equ(2.2)}) and (\ref{Equ(2.5)}), with some $c_1>0$ and $c_2>0$, we have  $\phi_\varepsilon(v_\varepsilon)\leq c_1$ and $|\phi_\varepsilon'(v_\varepsilon)|\leq c_2 v^{\alpha-1}_\varepsilon$ in $\Omega\times(0,\infty)$  for all $\varepsilon\in(0,1)$.
Then we utilize (\ref{Equ(2.5)}), $\alpha\geq\frac{1}{2}$ and Young's inequality to obtain
\begin{equation*}
\begin{aligned}
\int_{\Omega}\left|\nabla\left(u_{\varepsilon} \phi_{\varepsilon}\left(v_{\varepsilon}\right)\right)\right|^{\frac{4}{3}}\leq &2^\frac{1}{3}\int_{\Omega} \phi_{\varepsilon}^{\frac{4}{3}}\left(v_{\varepsilon}\right)\left|\nabla u_{\varepsilon}\right|^{\frac{4}{3}}+2^\frac{1}{3}c_2^\frac{4}{3}\int_{\Omega} u_{\varepsilon}^{\frac{4}{3}} v_{\varepsilon}^{\frac{4(\alpha-1)}{3}}\left|\nabla v_{\varepsilon}\right|^{\frac{4}{3}} \\
\leq&2^\frac{1}{3}\int_{\Omega} \phi_{\varepsilon}\left(v_{\varepsilon}\right) \frac{\left|\nabla u_{\varepsilon}\right|^2}{u_{\varepsilon}}+2^\frac{1}{3}\int_{\Omega} u_{\varepsilon}^2 \phi_{\varepsilon}^2\left(v_{\varepsilon}\right)\\
&+2^\frac{1}{3}c_2^\frac{4}{3}\left(\int_{\Omega}\frac{|\nabla v_{\varepsilon}|^4}{v_{\varepsilon}^3}+ \int_{\Omega} u_{\varepsilon}^2 v_{\varepsilon}^{\frac{4\alpha-1}{2}}\right) \\
\leq&2^\frac{1}{3}\int_{\Omega} \phi_{\varepsilon}\left(v_{\varepsilon}\right) \frac{\left|\nabla u_{\varepsilon}\right|^2}{u_{\varepsilon}}+2^\frac{1}{3}c_1^2\int_{\Omega} u_{\varepsilon}^2 \\
&+2^\frac{1}{3}c_2^\frac{4}{3}\left(\int_{\Omega}\frac{|\nabla v_{\varepsilon}|^4}{v_{\varepsilon}^3}+\|v_0\|^{\frac{4\alpha-1}{2}}_{L^\infty(\Omega)}\int_{\Omega} u_{\varepsilon}^2\right) \quad \text{for all } t>0 \text{ and } \varepsilon \in(0,1).
\end{aligned}
\end{equation*}
This together with (\ref{Equ(2.4)}), (\ref{Equ(2.15)}) and (\ref{Equ(4.22)}) yields the claimed assertion.  $\hfill{} \Box$

Collecting Corollary \ref{result4.02}, Lemmata \ref{result4.7}-\ref{result4.9}, \ref{result2.2} and \ref{result2.6}, we can extract suitably convergent subsequences of $(u_\varepsilon, v_\varepsilon)$ and a corresponding limit pair $(u, v)$ which solves (\ref{Equ(1.2)}) in the intended sense.
\begin{lemma}\label{result4.9}
Let $\gamma=2$, and assume (\ref{Equ(1.3)}), (\ref{Equ(1.5)}), (\ref{Equ(1.6)}) and (\ref{Equ(1.12)}) with $\alpha\geq\frac{1}{2}$. Then there exist $\left(\varepsilon_j\right)_{j \in \mathbb{N}} \subset(0,1)$ as well as functions $u$ and $v$ which satisfy (\ref{Equ(1.14)}) with $u \geq 0$ and $v>0$ a.e. in $\Omega \times(0, \infty)$, such that $\varepsilon_j \searrow 0$ as $j \rightarrow \infty$ and that beyond (\ref{Equ(4.13)})-(\ref{Equ(4.17)}) from Lemma \ref{result4.6} we have
\begin{equation*}
\nabla(\phi(v_{\varepsilon})u_{\varepsilon}) \rightharpoonup \nabla(\phi(v)u) \quad \text { in } L^{\frac{4}{3}}_{\text {loc }}(\bar{\Omega} \times[0, \infty))
\end{equation*}
as $\varepsilon=\varepsilon_j \searrow 0$, and that $(u, v)$ is a global weak solution of (\ref{Equ(1.2)}) according to Definition \ref{result0.1}.
\end{lemma}
\noindent{\bf{Proof.}} According to Lemmata \ref{result4.7}-\ref{result4.9}, \ref{result2.6}, \ref{result2.2} and Corollary \ref{result4.02}, we take a similar procedure as the proof of Lemma \ref{result4.6} to complete this proof. $\hfill{} \Box$

We can now easily prove Theorem \ref{result1.3}.\\
{\bf Proof of Theorem \ref{result1.3}.} Theorem \ref{result1.3} is asserted by the statement from Lemma \ref{result4.9}.
$\hfill{} \Box$

\section{Eventual smoothness}
\hspace*{\parindent} In this section, we will construct a Lyapunov functional to establish the eventual smoothness.  To achieve our goals, we follow an idea already introduced in \cite{LM3,W} and give the following two lemmas which are crucial to derive the lower bound of the cell mass.
\begin{lemma}\label{result5.1}
Let (\ref{Equ(1.3)}) hold and suppose $p>1$. Then we see that
\begin{equation}\label{Equ(5.1)}
\int_{0}^{\infty}\int_\Omega v^{p-2}_\varepsilon|\nabla v_\varepsilon|^2\leq\frac{1}{p(p-1)}\int_\Omega v^{p}_0 \quad \text{ for all } \;\; \varepsilon\in(0,1).
\end{equation}
\end{lemma}
\noindent{\bf{Proof.}} Multiplying the second equation of (\ref{Equ(2.1)}) by $v_\varepsilon^{p-1}$ and integrating over $\Omega$,  we have
\begin{equation*}
\begin{aligned}
\frac{d}{dt}\int_{\Omega} v_{\varepsilon}^{p}=-p(p-1)\int_{\Omega} v_{\varepsilon}^{p-2}|\nabla v_\varepsilon|^2-p\int_{\Omega}\frac{u_{\varepsilon}v_{\varepsilon}^{p}}{1+\varepsilon u_{\varepsilon}}
\end{aligned}
\end{equation*}
for all $t>0$ and $\varepsilon\in(0,1)$, which implies
\begin{equation*}
\begin{aligned}
\int_{\Omega} v_{\varepsilon}^{p}(\cdot,t)+p(p-1)\int_{0}^{t}\int_{\Omega} v_{\varepsilon}^{p-2}|\nabla v_\varepsilon|^2\leq\int_{\Omega} v_0^{p}
\end{aligned}
\end{equation*}
for all $t>0$ and $\varepsilon\in(0,1)$, this immediately shows (\ref{Equ(5.1)}).  $\hfill{} \Box$

\begin{lemma}\label{result5.2}
Let $\gamma=2$, and suppose that (\ref{Equ(1.3)}), (\ref{Equ(1.5)}) and (\ref{Equ(1.6)}) with $\alpha>0$ hold. Let
\begin{equation}\label{Equ(5.2)}
\mathcal{E_\varepsilon}(t):= \int_{\Omega}\left(u_\varepsilon(\cdot,t)-\frac{a}{b}-\frac{a}{b}\ln \frac{bu_\varepsilon(\cdot,t)}{a}\right) \quad \text{for all} \;\;t\geq0 \text{ and } \varepsilon\in(0,1).
\end{equation}
Then we obtain that $\mathcal{E_\varepsilon}(t)\geq0$ for all $t\geq0$ and $\varepsilon\in(0,1)$ and that there is $C_1>0$ fulfilling
\begin{equation}\label{Equ(5.3)}
\frac{d}{dt}\mathcal{E_\varepsilon}(t)+b \int_{\Omega}\left(u_\varepsilon-\frac{a}{b}\right)^2\leq C_1\int_{\Omega}v_\varepsilon^{\alpha-2}|\nabla v_\varepsilon|^2\quad \text{for all} \;\;t>0 \text{ and } \varepsilon\in(0,1).
\end{equation}
\end{lemma}
\noindent{\bf{Proof.}} We utilize Taylor's formula to see that $\mathcal{E_\varepsilon}(t)\geq0$ for all $t\geq0$ and $\varepsilon\in(0,1)$ (see \cite[Lemma 3.2]{BW} for details). Based on (\ref{Equ(1.5)}), (\ref{Equ(1.6)}), (\ref{Equ(2.2)}) and (\ref{Equ(2.5)}), we can fix positive constants $c_1$  and $c_2$  such that
$\phi_\varepsilon(v_\varepsilon)\geq c_1v^\alpha_\varepsilon$ and $ |\phi_\varepsilon'(v_\varepsilon)|\leq c_{2} v^{\alpha-1}_\varepsilon$ in $ \Omega\times(0,\infty)$   for all  $\varepsilon\in(0,1)$. By a direct calculation we infer that
\begin{equation}\label{Equ(5.4)}
\begin{aligned}
\frac{d}{dt}\mathcal{E_\varepsilon}(t)
=&a\int_{\Omega}u_\varepsilon-b\int_{\Omega}u_\varepsilon^2-\frac{a^2|\Omega|}{b}+a\int_{\Omega}u_\varepsilon-\frac{a}{b}\int_{\Omega}\phi_\varepsilon(v_\varepsilon)\frac{|\nabla u_\varepsilon|^2}{u_\varepsilon^2}\\
&-\frac{a}{b}\int_{\Omega}\phi'_\varepsilon(v_\varepsilon)\frac{\nabla u_\varepsilon\cdot\nabla v_\varepsilon}{u_\varepsilon}\\
\leq&-b \int_{\Omega}\left(u_\varepsilon-\frac{a}{b}\right)^2-\frac{a}{2b}\int_{\Omega}\phi_\varepsilon(v_\varepsilon)\frac{|\nabla u_\varepsilon|^2}{u_\varepsilon^2}+\frac{a}{2b}\int_{\Omega}\frac{|\phi'_\varepsilon(v_\varepsilon)|^2}{\phi_\varepsilon(v_\varepsilon)}|\nabla v_\varepsilon|^2\\
\leq&-b \int_{\Omega}\left(u_\varepsilon-\frac{a}{b}\right)^2+\frac{ac_2^2}{2bc_1}\int_{\Omega}v_\varepsilon^{\alpha-2}|\nabla v_\varepsilon|^2
\end{aligned}
\end{equation}
for all $t>0$ and $\varepsilon\in(0,1)$, which implies (\ref{Equ(5.3)}).  $\hfill{} \Box$

With Lemmata \ref{result5.1} and \ref{result5.2} at hand, we can derive the uniform lower bound of the cell mass.
\begin{lemma}\label{result5.3}
Let $\gamma=2$, and let (\ref{Equ(1.3)}), (\ref{Equ(1.5)}) and (\ref{Equ(1.6)}) with $\alpha>1$ be satisfied. Then there exists $C>0$ such that
\begin{equation}\label{Equ(5.5)}
\begin{aligned}
\int_{\Omega} u_\varepsilon(\cdot,t) \geq\frac{a|\Omega|}{be} \quad \text { for all }\,\,t>0 \text { and } \varepsilon\in(0,1)
\end{aligned}
\end{equation}
and
\begin{equation}\label{Equ(5.6)}
\int_{0}^{\infty}\int_{\Omega}\left(u_\varepsilon-\frac{a}{b}\right)^2\leq C\quad \text{for all} \;\; \varepsilon\in(0,1).
\end{equation}
\end{lemma}
\noindent{\bf{Proof.}} We apply (\ref{Equ(5.2)}), (\ref{Equ(5.3)}) and Lemma \ref{result5.1} to $p:=\alpha$ to find $c_1>0$ such that
\begin{equation}\label{Equ(5.7)}
\begin{aligned}
&\int_{\Omega}\left(u_\varepsilon-\frac{a}{b}-\frac{a}{b}\ln \frac{bu_\varepsilon}{a}\right)+b\int_{0}^{t}\int_{\Omega}\left(u_\varepsilon-\frac{a}{b}\right)^2\\
\leq&\int_{\Omega}\left(u_0-\frac{a}{b}-\frac{a}{b}\ln \frac{bu_0}{a}\right)+c_1\int_{0}^{t}\int_{\Omega}v_\varepsilon^{\alpha-2}|\nabla v_\varepsilon|^2\\
\leq&\int_{\Omega}\left(u_0-\frac{a}{b}-\frac{a}{b}\ln \frac{bu_0}{a}\right)+\frac{c_1}{\alpha(\alpha-1)}\int_\Omega v^{\alpha}_0
\end{aligned}
\end{equation}
for all $t>0$ and $\varepsilon\in(0,1)$,  along with the nonnegativity of $\mathcal{E_\varepsilon}(t)$ for all $t\geq0$ and $\varepsilon\in(0,1)$, we can see that (\ref{Equ(5.6)}) is valid. An application of (\ref{Equ(5.7)}) entails
\begin{equation*}
\begin{aligned}
\frac{a}{b}\int_{\Omega}\ln \frac{bu_\varepsilon}{a}\geq&\int_{\Omega}u_\varepsilon-\frac{a|\Omega|}{b}+b\int_{0}^{t}\int_{\Omega}\left(u_\varepsilon-\frac{a}{b}\right)^2\geq-\frac{a|\Omega|}{b}
\end{aligned}
\end{equation*}
for all $t>0$ and $\varepsilon\in(0,1)$, which implies
\begin{equation}\label{Equ(5.8)}
\begin{aligned}
\int_{\Omega}\ln u_\varepsilon\geq\left(\ln \frac{a}{b}-1\right)|\Omega|=\int_{\Omega}\ln \frac{a}{be}
\end{aligned}
\end{equation}
for all $t>0$ and $\varepsilon\in(0,1)$. Hence,
\begin{equation*}
\begin{aligned}
\int_{\Omega}u_\varepsilon\geq \frac{a|\Omega|}{be}
\end{aligned}
\end{equation*}
for all $t>0$ and $\varepsilon\in(0,1)$, this immediately gives (\ref{Equ(5.5)}).$\hfill{} \Box$

With the aid of (\ref{Equ(5.5)}), we can derive the decay of $\|v_\varepsilon(\cdot,t)\|_{L^\infty (\Omega)}$ by semigroup argument.
\begin{lemma}\label{result5.4}
Let $\gamma=2$, and assume that (\ref{Equ(1.3)}), (\ref{Equ(1.5)}) and (\ref{Equ(1.6)}) with $\alpha>1$ hold. Let $\left(\varepsilon_{j}\right)_{j \in \mathbb{N}}$ be as provided by Lemma \ref{result4.9}. Then whenever $\delta>0$, there exist $T(\delta)>0$ and $\varepsilon_0(\delta)\in(0,1)$ such that
\begin{equation*}
\begin{aligned}
\|v_\varepsilon(\cdot,t)\|_{L^\infty (\Omega)} \leq\delta \quad \text { for all }\,\,t>T(\delta) \text { and } \varepsilon\in(\varepsilon_{j})_{j \in \mathbb{N}} \cap(0,\varepsilon_0(\delta)).
\end{aligned}
\end{equation*}
\end{lemma}
\noindent{\bf{Proof.}} This assertion can be proved by an idea in the proof of \cite[Lemmata 5.3-5.5 and Corollary 5.1]{W}.$\hfill{} \Box$

Because of the possible diffusion degeneracy in the first equation from (\ref{Equ(2.1)}), it seems that the study of a differential inequality for $\int_{\Omega} \frac{u^{p}_{\varepsilon}(\cdot,t)}{(2 \delta-v_{\varepsilon}(\cdot,t))^{\theta}}$, functional of this type has successfully been utilized in \cite{LWR,TWE,WS}, is no longer available here. Fortunately, with the decay of $v_\varepsilon$ at hand, we can derive $L^p$-bounds of $u_\varepsilon$ for all $p>1$ after some waiting time by the construction of $\int_{\Omega}u^p_\varepsilon(\cdot,t)+\int_{\Omega} v^{-2p+1}_{\varepsilon}(\cdot,t)\left|\nabla v_{\varepsilon}(\cdot,t)\right|^{2p}$.
\begin{lemma}\label{result5.5}
Let $\gamma=2$, and let (\ref{Equ(1.3)}), (\ref{Equ(1.4)}) and (\ref{Equ(1.6)}) with $\alpha>1$ hold. Suppose that $\left(\varepsilon_{j}\right)_{j \in \mathbb{N}}$ is as defined by Lemma \ref{result4.9}. Then for all $p>1$, there exist $T(p)>0$, $\varepsilon_0(p)\in(0,1)$ and $C(p)>0$ such that
\begin{equation*}
\begin{aligned}
\|u_\varepsilon(\cdot,t)\|_{L^p(\Omega)} \leq C(p)\quad \text {for all}\,\,t>T(p) \text { and } \varepsilon\in(\varepsilon_{j})_{j \in \mathbb{N}} \cap(0,\varepsilon_0(p)).
\end{aligned}
\end{equation*}
\end{lemma}
\noindent{\bf{Proof.}} Given $p>1$, we choose a small $\delta=\delta(p)>0$ such that
\begin{equation}\label{Equ(5.9)}
\begin{aligned}
\kappa_1(p,n)C_2^{\frac{p+1}{p}}
\delta^\frac{\alpha(p+1)-1}{p}+\kappa_2(p,n)\delta\leq b,
\end{aligned}
\end{equation}
where $C_2$ is given by Lemma \ref{result3.1}, and $\kappa_1(p,n)$ and $\kappa_2(p,n)$ are taken from (\ref{Equ(3.22)}) and (\ref{Equ(3.23)}), respectively. An application of Lemma \ref{result5.4} implies the existences of $T(p)>0$ and $\varepsilon_0(p)>0$ such that
\begin{equation}\label{Equ(5.10)}
\begin{aligned}
\|v_\varepsilon(\cdot,t)\|_{L^\infty (\Omega)} \leq\delta \quad \text {for all}\,\,t>T(p) \text { and } \varepsilon\in(\varepsilon_{j})_{j \in \mathbb{N}} \text { such that } \varepsilon<\varepsilon_0(p).
\end{aligned}
\end{equation}
Due to $\gamma=2$, it follows from Lemma \ref{result3.2} with the choice $q:=2p$, (\ref{Equ(3.22)}) and (\ref{Equ(3.23)}) that one can find some $c>0$ such that
\begin{equation}\label{Equ(5.11)}
\begin{aligned}
\frac{d}{dt}&\left(\int_{\Omega}u_\varepsilon^p+\int_{\Omega} v_{\varepsilon}^{-2p+1}\left|\nabla v_{\varepsilon}\right|^{2p}\right)+\int_{\Omega}u_\varepsilon^p+\int_{\Omega} v_{\varepsilon}^{-2p+1}\left|\nabla v_{\varepsilon}\right|^{2p}\\
\leq& \frac{p(p-1)}{2}\left(\frac{2p-1}{(p-1)(2p+\sqrt{n})^2}\right)^{-\frac{1}{p}}C_2^{\frac{p+1}{p}}
\|v_\varepsilon(\cdot,t)\|_{L^\infty(\Omega)}^\frac{\alpha(p+1)-1}{p}\int_{\Omega}u_\varepsilon^{p+1}\\
&+\frac{2^{3p+1}p(2p+\sqrt{n}+1)^{2p}(2p-2+\sqrt{n})^{p+1}\|v_\varepsilon(\cdot,t)\|_{L^\infty(\Omega)}}{(2p-1)^p}
\int_{\Omega}u_\varepsilon^{p+1}-\frac{bp}{2}\int_{\Omega} u_\varepsilon^{p+1}+c\\
=& -\frac{p}{2}\left\{b-\kappa_1(p,n)C_2^{\frac{p+1}{p}}
\|v_\varepsilon(\cdot,t)\|_{L^\infty(\Omega)}^\frac{\alpha(p+1)-1}{p}-\kappa_2(p,n)\|v_\varepsilon(\cdot,t)\|_{L^\infty(\Omega)}\right\}\int_{\Omega}u_\varepsilon^{p+1}+c
\end{aligned}
\end{equation}
for all $t>0$ and $\varepsilon\in(0,1)$.  Hence, we utilize (\ref{Equ(5.9)})-(\ref{Equ(5.11)}) and an ODE argument to obtain that  $\|u_\varepsilon(\cdot,t)\|_{L^p(\Omega)}$ is bounded for all $t>T(p)$ and $\varepsilon\in(\varepsilon_{j})_{j \in \mathbb{N}}$ such that $\varepsilon<\varepsilon_0(p).$ Then the proof of Lemma \ref{result5.5} is completed. $\hfill{} \Box$

With the $L^p$-estimates of $u_\varepsilon$ for arbitrary large $p>1$ at hand, the higher order regularity features of $u_\varepsilon$ and $v_\varepsilon$ can be derived.
\begin{lemma}\label{result5.6}
Let $\gamma=2$, and suppose that (\ref{Equ(1.3)}) and (\ref{Equ(1.4)})-(\ref{Equ(1.6)}) with $\alpha>1$ hold. Let $\left(\varepsilon_{j}\right)_{j \in \mathbb{N}}$ be as given by Lemma \ref{result4.9}. Suppose that $(u_\varepsilon,v_\varepsilon)$ is a global solution of (\ref{Equ(2.1)}). Then there exists $T_0>0$ such that for any $T>T_0$, there are $\sigma=\sigma(T) \in(0,1)$, $\varepsilon_0(T)>0$ and $C(T)>0$ such that
\begin{equation*}\label{Equ(7.24)}
\|u_\varepsilon\|_{C^{2+\sigma, 1+\frac{\sigma}{2}}(\bar{\Omega} \times[T_0, T])}+\|v_\varepsilon\|_{C^{2+\sigma, 1+\frac{\sigma}{2}}(\bar{\Omega} \times[T_0, T])} \leq C(T)
\end{equation*}
for all $\varepsilon\in(\varepsilon_{j})_{j \in \mathbb{N}}$ such that $\varepsilon<\varepsilon_0(T).$ Moreover, there exists a subsequence $\left(\varepsilon_{j_k}\right)_{k \in \mathbb{N}}$ of $\left(\varepsilon_{j}\right)_{j \in \mathbb{N}}$ such that
\begin{equation*}\label{Equ(7.25)}
u_\varepsilon\rightarrow u,\,\,v_\varepsilon\rightarrow v \text { in } C_{loc}^{2, 1}(\bar{\Omega} \times[T_0, \infty)) \,\,\text{as } \varepsilon=\varepsilon_{j_k}\searrow 0.
\end{equation*}
\end{lemma}
\noindent{\bf{Proof.}} This proof can be proved by taking almost procedure as in Lemmata \ref{result3.9}, \ref{result3.11} and \ref{result3.12}
(see also \cite[Lemma 4.5]{WAA} or \cite[Lemmas 4.4-4.6]{LWR}).  $\hfill{} \Box$

Now, we begin to prove Theorem \ref{result1.4}. \\
{\bf Proof of Theorem \ref{result1.4}.} Theorem \ref{result1.4} immediately results from Lemma \ref{result5.6}.\\
{\bf Acknowledgments}

The author warmly thanks Professor Michael Winkler for his carefully reading and valuable suggestions which greatly improved this work.
This work is supported by Natural Science Foundation of Chongqing (No. cstc2021jcyj-msxmX0412) and the NNSF of China (No. 12271064).

\end{document}